\documentclass[12pt]{amsart}
\usepackage{latexsym}
\usepackage{amssymb}
\usepackage{amsxtra}
\usepackage{amsmath}
\usepackage{amsfonts}

\newcommand{\beq}{\begin{equation}}
\newcommand{\eeq}{\end{equation}}
\newcommand{\beqa}{\begin{eqnarray}}
\newcommand{\eeqa}{\end{eqnarray}}
\newcommand{\beaa}{\begin{eqnarray*}}
\newcommand{\ben}{\begin{eqnarray*}}
\newcommand{\eaa}{\end{eqnarray*}}
\newcommand{\een}{\end{eqnarray*}}

\newcommand \nc {\newcommand}

\newtheorem{theorem}{Theorem}[section]
\newtheorem{lemma}[theorem]{Lemma}
\newtheorem{proposition}[theorem]{Proposition}
\newtheorem{corollary}[theorem]{Corollary}
\newtheorem{definition}[theorem]{Definition}

\newtheorem{remark}[theorem]{Remark}
\newtheorem{conjecture}[theorem]{Conjecture}

\nc \thref{Theorem \ref}
\nc \leref{Lemma \ref}
\nc \prref{Proposition \ref}
\nc \coref{Corollary \ref}
\nc \deref{Definition \ref}
\nc \exref{Example \ref}
\nc \reref{Remark \ref}

\newcommand{\CC}{\mathcal{C}}
\newcommand{\X}{\mathcal{X}}

\newcommand{\W}{\mathcal{W}}
\newcommand{\A}{\mathcal{A}}
\newcommand{\B}{\mathcal{B}}
\newcommand{\C}{\mathbb{C}}
\newcommand{\D}{\mathcal{D}}

\renewcommand{\H}{\mathcal{H}}
\newcommand{\J}{\mathcal{J}}
\renewcommand{\L}{\mathcal{L}}
\newcommand{\M}{\mathcal{M}}

\newcommand{\QQ}{\mathbb{Q}}

\newcommand{\T}{\mathcal{T}}

\newcommand{\Z}{\mathbb{Z}}

\newcommand{\f}{\mathbf{f}}

\newcommand{\q}{\mathbf{q}}
\renewcommand{\t}{\mathbf{t}}
\newcommand{\x}{\mathbf{x}}
\newcommand{\y}{\mathbf{y}}

\newcommand{\Mbar}{\overline{\mathcal{M}}}

\def\res{\mathop{\rm Res}\nolimits}

\def\d{\partial}

\def\iso{\cong}
\def\tensor{\otimes}

\def\Poincare{Poincar\'e}

\def\({\left(}
\def\){\right)}
\def\[{\left[}
\def\]{\right]}
\def\<{\left\langle}
\def\>{\right\rangle}

\def\liegl{{\mathfrak{gl}}}

\def\lieGL{{\rm GL}}

\def\gl{\lambda}
\def\ge{\epsilon}
\def\ga{\alpha}
\def\gd{\delta}
\def\gb{\beta}

\newcommand{\bpsi}{\bar{\psi}}
\newcommand{\proj}{\mathbb{P}}

\newcommand{\com}{\mathbb{C}}

\newcommand{\mC}{\mathcal{C}}

\newcommand{\sI}{\mathcal{I}}

\newcommand{\bt}{\mathbf{t}}
\newcommand{\bq}{\mathbf{q}}

\newcommand{\sO}{\mathcal{O}}

\newcommand\radice[2]{\sqrt[\uproot{2}#1]{#2}}

\title[Laurent polynomials, $\mathbb{P}^1$-orbifolds, and integrable hierarchies]
{The spaces of Laurent polynomials, Gromov-Witten theory of $\mathbb{P}^1$-orbifolds, and 
integrable hierarchies}

\author{Todor E. Milanov}
\address{Department of Mathematics\\ Stanford University\\ 
Stanford\\ CA 94305--2125\\ USA}
\email{milanov@math.stanford.edu}

\author{Hsian-Hua Tseng}
\address{Department of Mathematics\\ University of British Columbia\\ 
1984 Mathematics Road\\ Vancouver\\ B.C. V6T 1Z2\\ Canada}
\curraddr{Department of Mathematics\\ University of Wisconsin-Madison\\ Van Vleck Hall, 480 Lincoln Drive\\ Madison, WI 53706-1388\\ USA}
\email{tseng@math.wisc.edu}

\date{\today}
\thanks{{\em 2000 Math. Subj. Class.} 14N35, 17B69, 32S30}

\begin{document}
\begin{abstract}
Let $M_{k,m}$ be the space of Laurent polynomials in one variable 
$x^k + t_1 x^{k-1}+
\ldots t_{k+m}x^{-m},$ where $k,m\geq 1$ are fixed integers and 
$t_{k+m}\neq 0.$
According to B. Dubrovin \cite{D}, $M_{k,m}$ can be equipped with a 
semi-simple Frobenius structure. 
In this paper we prove that the corresponding descendent  and ancestor 
potentials of $M_{k,m}$ (defined as in \cite{G1}) satisfy Hirota 
quadratic equations (HQE for short). 

Let $\CC_{k,m}$ be the orbifold obtained from $\mathbb{P}^1$ by cutting
small discs $D_1\iso \{|z|\leq \ge\}$ and $D_2\iso\{|z^{-1}|\leq \ge\}$ around 
$z=0$ and $z=\infty$ and 
gluing back the orbifolds $D_1/\Z_k$ and $D_2/\Z_m$ in the obvious way.
We show that the orbifold quantum cohomology of $\CC_{k,m}$ coincides 
with $M_{k,m}$ as Frobenius manifolds. Modulo some yet-to-be-clarified 
details, this implies that the descendent 
(respectively the ancestor) potential of $M_{k,m}$ is a generating 
function for the descendent (respectively ancestor) orbifold Gromov--Witten 
invariants of $\CC_{k,m}$.

There is a certain similarity between
our HQE and the Lax operators of the Extended bi-graded Toda hierarchy, 
introduced by G. Carlet in \cite{car}. 
Therefore, it is plausible that our HQE characterize the tau-functions 
of this hierarchy and we expect that the Extended bi-graded Toda hierarchy
governs the Gromov--Witten theory of $\CC_{k,m}.$ 

\end{abstract}
\maketitle


\section{Introduction}

\subsection{Background}

By definition (see \cite{D} or \cite{Ma}), a Frobenius structure on a manifold $M$ is
a collection of a flat metric $g$ on $M,$ a multiplication $\bullet$ in each 
tangent space $T_tM,$ depending smoothly on $t$ and satisfying the Frobenius 
property 
$g(X\bullet Y,Z)=g(Y,X\bullet Z),$ and a flat vector field $e$ which is a unity 
with respect to $\bullet,$ such that certain integrability conditions are
satisfied. For example,
if $X$ is a compact symplectic manifold then the cohomology algebra 
$H^*(X)$ is naturally equipped with a Frobenius structure where the
metric is given by the {\Poincare} pairing and the multiplication by the quantum cup product, see e.g. \cite{G3} for more details. 

A Frobenius manifold $M$ is called (generically) {\em semi-simple} if there
exists a point $t\in M$ such that the corresponding tangent space $T_tM$ is a 
semi-simple algebra, i.e., it has no nilpotents. For semi-simple $M$,
A. Givental \cite{G1} introduced the so-called {\em total descendent} and 
{\em total ancestor} potentials, denoted respectively by $\D^M$ and 
$\A_t^M,$ where $t\in M$ is a semi-simple point. They belong to the {\em Fock space} $B$, which is an infinite dimensional vector space described as 
follows: if we pick a trivialization of the tangent bundle $TM\iso M\times H,$ 
corresponding to a choice of a flat coordinate system on $M$ (here $H$ is 
an arbitrary fixed tangent space of $M$), then $B$ is a certain 
completion of the space of functions on $\H_+:=H[z].$ Moreover, A. Givental 
conjectured that if $M$ is a Frobenius manifold coming from the quantum cohomology theory of
a compact K\"ahler manifold $X$ then $\D^M$ (respectively $\A_t^M$) are 
generating functions for the descendent (respectively ancestor) Gromov--Witten invariants of $X$. This conjecture is proven for toric manifolds (see \cite{G3} for Toric Fano case and \cite{Ir} for general toric case), Flag manifolds \cite{KJ}, and Grassmannians \cite{BCFK}. Recently, C. Teleman announced a classification of semi-simple cohomological field theories. Together with some yet-to-be-clarified technical details, this implies Givental's conjecture in general. 

Let $q_n^a,$ $1\leq a\leq N$, $n=0,1,2,\ldots$ be a set of formal variables. 
By fixing a basis in $\H_+,$ we identify the Fock space $B$ with the 
space of formal series on the variables $q_n^a$ with complex coefficients. 
Given $\tau\in B,$ we will refer to the coefficients in the corresponding formal series
as {\em Fourier coefficients}. In this paper we prove that if $M=M_{k,m}$ is
the space of Laurent polynomials in one variable, then  the Fourier 
coefficients of $\D^M$ and $\A^M_t$ satisfy an infinite system of quadratic
relations. Alternatively, these quadratic relations can be written as an 
infinite system of PDE's which involve quadratic expressions of $\tau$, its
partial derivatives, and its translation. We refer to such 
a system of PDEs as Hirota Quadratic Equations (HQE for short). Recently, G. Carlet \cite{car} associated an integrable hierarchy to $M_{k,m}$ which fits in the general 
framework of \cite{DZ2}. We expect that our HQE give a description of
Carlet's hierarchies in terms of HQEs and tau-functions.
We also prove that the Frobenius manifold $M_{k,m}$ is isomorphic to the 
orbifold quantum cohomology of $\CC_{k,m}$.

Frobenius manifolds and integrable systems are closely related (c.f. \cite{DZ2}). Some classes of integrable systems can be described in terms of {\em $\tau$-functions} and Hirota quadratic 
equations (also known as Hirota bi-linear equations). Examples include KdV, KP, and Toda lattice hierarchies.
Here is, to the best of our knowledge,  a complete list of pairs 
consisting of a semi-simple Frobenius manifold $M$ and an integrable hierarchy for
which it is known that the potential $\D^M$ is a tau-function of the corresponding
hierarchy:
\begin{itemize}
\renewcommand{\labelitemi}{---}
\item
$M$ is the space of miniversal deformations of $A,$ $D,$ or $E$ type singularity 
AND
the Kac--Wakimoto hierarchies corresponding to the Coxeter transformation
in the Weyl algebra of the simple Lie algebras of $A,$ $D,$ or $E$ type (\cite{G2}, \cite{GM}).
\item
Quantum cohomology of a point (which coincides with the miniversal
deformation of $A_1$ singularity) AND the KdV hierarchy (\cite{W}, \cite{Ko}).
\item
Quantum cohomology of $\C \proj^1$ AND the Extended Toda hierarchy (\cite{M}, \cite{OP1}).

\item
Equivariant quantum cohomology of $\C \proj^1$ AND the 2-Toda hierarchy (\cite{M3}, \cite{OP1}).
\item
Orbifold quantum cohomology of the classifying stack $BG$ of a finite group $G$
AND $|Conj(G)|$ commuting copies of the KdV hierarchies, where $Conj(G)$ is the set of conjugacy classes of $G$ (\cite{JK}).
\end{itemize}
The results of this paper suggest that we can add one more pair to the above list:
\begin{itemize}
\renewcommand{\labelitemi}{---}
\item
Orbifold quantum cohomology of $\CC_{k,m}$ AND the Extended bi-graded Toda
hierarchy.  
\end{itemize}

To complete this, it remains to clarify the following details: the functions satisfying 
our HQEs 
are tau-functions of the Extended bi-graded Toda hierarchy, and the potential 
$\D^M$ (respectively $\A_t^M$) is a generating function for descendent (respectively ancestor) orbifold Gromov--Witten invariants of $\CC_{k,m}$. The solutions to these two problems should not be very difficult: for the first one
we need to generalize the techniques from \cite{M2}, and the second one
follows either from Teleman's work \cite{te}, or alternatively can be proven by virtual localization (see \cite{ts} for details).

\subsection{Summary of results}
Given positive integers $k$ and $m$ and a non-zero complex number $Q$, 
we denote by $M$ the space of Laurent polynomials 
\ben
f = x^k+\sum_{i=1}^k t_ix^{k-i} + \sum_{j=1}^{m-1} t_{k+j} \(Qe^{t_N}/x\)^j 
+\(Qe^{t_N}/x\)^m,
\een
where $N=k+m,$ i.e., $M\iso \C^{N-1}\times \C^*.$ 
Each tangent space $T_fM$ is naturally identified with the local 
algebra $\C [x,x^{-1}]/\<\d_xf\>$: the vector field 
$\d/\d t_i$ corresponds to the projection of $\d f/\d t_i$ in 
$\C [x,x^{-1}]/\<\d_xf\>.$ Via this identification the product in the 
local algebra defines an associative, commutative product $\bullet_f$
on the tangent space $T_fM$ with unity $e=\d/\d{t_k}$. 

Furthermore, let $\omega = dx/x$ be the standard volume form on $\C^*.$ 
Then we equip each tangent space with a residue pairing
\beq\label{res_metric}
\(\d_{t_i},\d_{t_j}\)_f =
-\( {\rm res}_{x=0} + {\rm res}_{x=\infty}\)
\frac{\( \d_{t_i}f\omega\)\( \d_{t_j}f\omega\)}{df}.
\eeq    
Finally, we assign degrees to $x$ and $t_i$ such that $f$ becomes 
a homogeneous polynomial of degree 1. In order to keep track of the 
homogeneity properties of functions on $M$, we introduce the
following {\em Euler vector field}: 
$$
E=\sum_{i=1}^{k} \frac{i}{k} t_i\d_{ t_i} + 
  \sum_{j=1}^{m-1}\(1-\frac{j}{m}\)t_{k+j}\d_{ t_{k+j}}+ 
  \(\frac{1}{k}+\frac{1}{m}\)\d_{ t_N}. 
$$

The data introduced here satisfy an integrability condition: for 
each $z\in \C^*$ 
\beq\label{flat_connection}
\nabla = \nabla^{\rm L.C.}-z^{-1}\sum_{i=1}^N (\d_{t_i}\bullet) d t_i, 
\eeq
is a flat connection (i.e. $\nabla^2 = 0$) on $TM,$ where $\nabla^{\rm L.C.}$
is the Levi-Civita connection of the residue metric.   
In particular, $\nabla^{\rm L.C.}$ is flat as well.

Let $A_{ij}^p$ be the structure constants of $\bullet$ and $g_{ij}$ is
the tensor of the residue metric, i.e.,
$\d/\d t_i\bullet_f \d/\d t_j = \sum_p A_{ij}^p(f)\d/\d t_p$ and 
$g_{ij}(f):=(\d/\d t_i,\d/\d t_j)_f.$ Note that  
$A_{ij}^p$ and $g_{ij}$ are polynomials in $t_1,\ldots ,t_{N-1},Qe^{t_N},$ 
thus taking the corresponding free terms yields (in each tangent 
space $T_fM$) an associative, commutative multiplication, which will be 
called {\em cup product} or {\em classical multiplication}, 
and a non-degenerate bilinear pairing.    
The corresponding algebra structure on $T_fM$ can be described 
explicitly as follows: under the map 
\ben
\d/\d t_{k-i} \mapsto  \phi_i:=X^i,\ \d/\d t_{k+j}\mapsto \phi_{k+j}:= Y^j ,
\d/\d t_N\mapsto m\phi_N:=mY^m,
\een 
where $1\leq i\leq k-1,\ 0\leq j\leq m-1,$
the cup product corresponds to the multiplication in the algebra
$H:= \C[X,Y]/ \< kX^k-mY^{m}, XY\>$ and the free terms of
$g_{ij}$ induce a non-degenerate bilinear pairing on $H:$ 
$(\phi_i ,\phi_{k-i} ) = 1/k,$ $(\phi_{k+j} ,\phi_{k+m-j} ) = 1/m,$
and all other pairs of vectors are orthogonal.  
Moreover, using the Levi--Civita connection we can choose a flat 
coordinate system $\tau_1,\ldots,\tau_{N-1},e^{\tau_N}$ on $M$ such that 
the map $\d/\d \tau_i\rightarrow \phi_i, 1\leq i\leq N,$ gives 
a trivialization of the tangent bundle $TM\iso M\times H$ under which 
the cup product and the residue pairing  correspond respectively 
to the multiplication and the bilinear pairing of $H$. Such flat coordinates 
will be constructed explicitly in Section 
\ref{sec:flat_structure}. We will denote by $\d_i$ the  vector 
field $\d/\d\tau_i$ and by $f_\tau$ the Laurent polynomial in $M$ 
corresponding to $\tau=(\tau_1,\ldots,\tau_N).$ 

A direct computation of the orbifold cohomology $H_{orb}^*(\CC_{k,m};\C)$
shows that 
$H_{orb}^*(\CC_{k,m};\C)\iso H$ where the isomorphism is given by  
$X={\rm P.D.}([B\mathbb{Z}_k])$ and $Y= {\rm P.D.}([B\mathbb{Z}_m]).$
Our first result is: 
\begin{theorem}\label{t3} 
 $M_{k,m}$ is isomorphic to the Frobenius manifold corresponding to the big orbifold quantum cohomology of $\CC_{k,m}$. 
\end{theorem} 
In other words $M_{k,m}$ is the full mirror model of $\CC_{k,m}.$ 

Details of the proof of Theorem \ref{t3} and some background on orbifolds and their quantum cohomologies will be given in Section \ref{sec:ocp1}. 

Let $\H = H((z^{-1}))$ be the space of formal Laurent series 
in $z^{-1}$ with vector coefficients equipped with a symplectic structure, 
\ben
\Omega(f,g) := \frac{1}{2\pi i} \oint \(f(-z),g(z)\)dz,\ f, g\in \H.
\een 
The polarization $\H=\H_+\oplus \H_-$, defined by the Lagrangian subspaces
$\H_+ = H[z]$ and $\H_- =z^{-1} H[[z^{-1}]]$, identifies $\H$ with the 
cotangent bundle $T^*\H_+$.

Let $\ge$ be a formal variable -- the genus parameter in Gromov--Witten theory.
By definition, {\em the (Bosonic) Fock space} $B_H$ is the vector space
of functions on $\H_+,$ completed in a certain way. Namely, if we let 
$\q(z)=\sum_{k\geq 0} q_kz^k\in \H_+$ then $B_H$ is the space of formal series
in the sequence of  vector variables $q_0,q_1+{\bf 1},q_2,\ldots ,$
whose coefficients are formal Laurent series in $\ge.$ 
We construct a representation of the {\em Heisenberg 
Lie algebra} generated by the linear Hamiltonians on the Fock space $B_H.$
Let $\{\phi^i\}$ be a basis of $H$ dual to $\{\phi_i\}$ with respect to
the residue pairing. Then the linear functions on $\H$ defined by 
$p_{k,i}=\Omega(\cdot, \phi_i z^k)$ and $q_{k}^{i}=\Omega(\phi^i(-z)^{-k-1},\cdot)$ 
form a Darboux coordinate system on $\H.$  Thus the formulas 
\beq
\label{qlinear}
\widehat{ q}_{k}^{i} := {q_{k}^{i}}/{\ge},\quad
\widehat{ p}_{k,i} := \ge {\d}/{\d q_{k}^{i}},
\eeq
define a representation on $B_H.$ 
Given a vector ${\bf f}\in \H$ we define 
{\em a vertex operator} acting on $B_H$: 
$e^{\widehat \f}:=(e^\f)\sphat:=e^{\hat \f_-}e^{\hat \f_+} ,$ where 
$\f_\pm$ is the projection of $\f$ on $\H_\pm$ and $\f_\pm$ is 
identified with the linear Hamiltonian $\Omega(\cdot ,\f_\pm ).$

A fundamental solution to the system of differential equations 
corresponding to the flat connection \eqref{flat_connection} has 
two singularities -- at $z=0$ and $z=\infty.$ The information 
about these singular points is encoded in two vectors  
$\D^M,\A_\tau^M\in B_H,$ called {\em total descendent} and 
{\em total ancestor} potentials of
$M,$  where $\tau$ is {\em a semi-simple point}, i.e., for $\tau'$ in a 
neighborhood
of $\tau$ the critical values of $f_{\tau'}$ form a local coordinate 
system on $M$. In other words, $f_\tau$ has only Morse type critical 
points, see Section \ref{sec:desc_anc} for precise definitions. 

\medskip

Let $\Delta\subset M\times \C$ be the set of pairs $(\tau,\lambda)$ such that 
the equation $f_\tau(x)=\lambda$ has less than $N$ solutions. Then the space
$\(M\times \C \)\backslash \Delta$ admits an $N$-fold covering  $V:$ the   
fiber over $(\tau,\lambda)$ is $V_{\tau,\lambda} := f_\tau^{-1}(\lambda).$
The relative homology groups $H_1(\C^*,V_{\tau,\gl};\C)$ vary naturally
with respect to $(\tau,\gl)$, so they define a vector bundle on 
$(M\times\C)\backslash \Delta.$ Moreover, this bundle is 
equipped with a connection, called {\em Gauss-Manin connection}: given a 
path $C$ from $(\tau_1,\gl_1)$ to $(\tau_2,\gl_2)$ there is a natural 
identification (since $V$ is a fibration) between the corresponding
relative homology groups. 

Fix an arbitrary reference point 
$(\tau_0,\gl_0)\in \(M\times \C \)\setminus \Delta.$   
For each $\gb\in H_1(\C^*,V_{\tau_0,\gl_0};\C)$ and $n\in \Z$ we define
multivalued period mappings $I_{\gb}^{(n)}:(M\times\C)\backslash \Delta\rightarrow H$ as follows: 
{\allowdisplaybreaks
\ben
(I^{(-p)}_{\gb}(\tau,\gl),\d_i) & = &
-\d_i\,\int_{\gb(\tau,\gl)} \frac{(\gl-f_\tau)^{p}}{p!}\omega, 
\quad 1\leq i\leq N \\
I^{(p)}_{\gb}(\tau,\gl)  & = & \d_\gl^{p}I_{\gb}^{(0)}(\tau,\gl),
\een}
where $p$ is a non-negative integer and 
$\gb(\tau,\gl)\in H_1(\C^*, f_\tau^{-1}(\gl);\Z)$ is a cycle obtained from
$\gb$ via a parallel transport along a path $C$ connecting $(\tau_0,\gl_0)$ 
and $(\tau,\gl).$ Note that the value of $I^{(n)}_\gb$ 
depends on the choice of the path $C.$ Finally, put
\ben
\f^\gb_\tau(\gl) = \sum_{n\in \mathbb{Z}} I^{(n)}_\gb(\tau,\gl)(-z)^n,\quad 
\Gamma_\tau^\gb = e^{\hat \f_\tau^\gb}.
\een

Let $x_1^0\,\ldots,x_N^0$ be the solutions to  $f_{\tau_0}(x)=\gl_0$.  For each 
$1\leq a\leq N$ choose a path in $\C^*$ from 1 to $x_a^0,$ i.e., fix a value 
of $\log x_a^0.$ We introduce vertex operators $\Gamma_\tau^a,$ 
$a=1,2,\ldots,N$ corresponding to the one-point cycles 
$[x_a^0]\in H_0(f_{\tau_0}^{-1}(\gl_0);\Z)$, as follows. Given an integer 
$n$, we define a multivalued period mapping 
$I_a^{(n)}:(M\times\C)\backslash \Delta \rightarrow H$ by 
\beqa\label{def:1pt}
&&
(I_a^{(-p)}(\tau,\gl),\d_i) = 
-\d_i \int_{[x_a]}
d^{-1}\( \frac{1}{p!} \(\gl-f_\tau\)^p \omega \), \\ \notag
&&
I^{(p)}_a = \d_\gl^p I^{(0)}_a,
\eeqa
where $p\geq 0$ and $d^{-1}$ is a linear operator acting on the space 
of volume forms on $\C^*$ according to the rule 
$d^{-1}(x^kdx) = x^{k+1}/(k+1)$ if $k\neq -1$ and $d^{-1}(dx/x)=\log x$.
The periods  $I_a^{(n)}$ are multi-valued: the values of $x_a=x_a(\tau,\gl)$ 
and $\log x_a$ depend on the choice of a path (avoiding $\Delta$) from 
$(\tau_0,\gl_0)$ to $(\tau,\gl).$ Finally, put
\ben
\f_\tau^a(\gl) = \sum_{n\in\Z} I_a^{(n)}(\tau,\gl)(-z)^n,\quad
\Gamma_\tau^a = e^{\hat\f_\tau^a}.
\een
Note that if $\gb\in H_1(\C^*,f_{\tau_0}^{-1}(\gl_0);\Z)$ is a relative
cycle represented by the composition of the two paths $x_a^0$ to 1 and 1 to 
$x_b^0$ -- the same ones which specify the branch of $\log x_a^0$ and
$\log x_b^0$, then a simple application of the Stokes' formula implies:
$\f^\gb_\tau = \f_\tau^b - \f_\tau^a.$

The vertex operators $\Gamma_\tau^a,\ 1\leq a\leq N$ depend on the 
choice of $\log x_a^0$ as follows. Let $\phi\in H_1(\C^*;(2\pi i)^{-1}\Z)$ be a
cycle normalized by $\int_\phi\omega =1.$ Then changing the value from
$\log x_a^0$ to  $\log x_a^0 + r_a,$ $r_a\in 2\pi i\Z,$ transforms the vertex 
operators 
\ben
\Gamma_\tau^a \tensor\Gamma_\tau^{-a} \quad\mbox{ into }\quad 
\(\Gamma_\tau^{r_a\phi} \tensor\Gamma_\tau^{-r_a\phi}\)
\(\Gamma_\tau^a \tensor\Gamma_\tau^{-a}\).
\een
To offset this ambiguity,
we allow vertex operators acting on a larger Fock space
$\B_H := \A\tensor_\C B_H$. Here $\A$ is the algebra of differential 
operators 
$\sum_{0\leq k \leq N} a_k(x;\ge)\d_x^k,$ where each
$a_k$ is a formal Laurent series in $\ge$ with coefficients smooth
functions in $x.$ We equip $\A$ with an {\em anti-involution}
$\#$ defined by its action on the generators $x$ and $\ge\d_x$ of 
$\A:$
\ben
(\ge\d_x)^\# = -\ge\d_x ,\quad x^\# = x.
\een
Let $w_\infty=-d\tau_k\,z^{-1}$ and $v_\infty =\d_k.$ There are unique
vectors $w_\tau,v_\tau\in \H$ such that they are horizontal sections
of $\nabla$ (see \eqref{flat_connection}), depend polynomially on 
$\tau_1,\ldots,\tau_N, Qe^{\tau_N}$ and their free terms are 
respectively $w_\infty$ and $v_\infty.$ Introduce a vertex 
operator (acting on $\B_H$)
\ben
\Gamma_\tau^\delta:= 
\exp\left( 
(\f_\tau^{\phi}-w_\tau)\ge\d_x \right)\sphat\ 
\exp \left( {x} v_\tau/\ge \right)\sphat.
\een
It has the following crucial property:
\beq \label{monodromy:killer}
\(\Gamma_\tau^{\delta \#}\tensor\Gamma_\tau^\delta\)
\(\Gamma_\tau^{r\,\phi}\tensor  \Gamma_\tau^{-r\,\phi}\) =
e^{(\hat w_\tau\tensor 1 - 1\tensor \hat w_\tau)\,r}
\Gamma_\tau^{\delta \#}\tensor\Gamma_\tau^\delta.
\eeq
Finally, for each $i$ with  $1\leq i\leq N$, define 
$c_\tau^i(\gl) = {1}/{f_\tau'(x_i)},$ where $x_i=x_i(\tau,\gl)$ is a 
solution to 
$f_\tau(x)=\gl$ and $'$ is the derivative with respect to $x.$ 
\begin{definition}{\em 
We say that $\T\in B_H$ satisfies the HQE 
\eqref{HQE:ancestors} if the 1-form
\beqa
\label{HQE:ancestors}
\(\Gamma_\tau^{\gd\#}\tensor \Gamma_\tau^\gd \) 
\(
\sum_{i=1}^N c_\tau^i \Gamma_\tau^{i}\tensor \Gamma_\tau^{-i} 
\)
\(\T\tensor \T\)\  {d\gl},
\eeqa
computed at $\q'$ and $\q''$ such that
$\widehat w_\tau' - \widehat w_\tau'' =r$, is regular in 
$\gl$ for each $r\in \Z$.}
\end{definition}
Here $\T\tensor\T$ means the function $\T(\q')\T(\q'')$ on the 
two copies of the variable $\q =\{ q_k^a\ |\ 1\leq a\leq N, k\geq 0\}$
and the vertex operators in $\Gamma_\tau^*\tensor\Gamma_\tau^{-*}$  preceding
(respectively following) $\tensor$ act on $\q'$ (respectively on $\q''$).
Furthermore, $\widehat w_\tau$ is identified with a linear function in $\q$
via the symplectic form, i.e., $\widehat w_\tau(\q) = \ge^{-1}\Omega(\q,w_\tau),$
where $\q:= \sum q_k^a\phi_az^k.$ 
The expression \eqref{HQE:ancestors} is interpreted as taking values 
in the vector space $\B_H\tensor_\A\B_H.$ Thanks to \eqref{monodromy:killer},
when $\widehat w_\tau'-\widehat w_\tau''=r\in \Z$ the expression
\eqref{HQE:ancestors} is single-valued near $\gl=\infty.$ After the change 
$\y=(\q'-{\q}'')/(2\ge), \x=(\q'+{\q}'')/2$ and the substitution\footnote{Note that 
$w_\tau = d\tau_k (-z)^{-1} + O(z^{-2}).$ Also 
$\widehat w_\tau(\y) = \widehat w_\tau(\q'-\q'')/(2\ge) = 
(\widehat w_\tau'-\widehat w_\tau'')/(2\ge) =r/(2\ge).$ Thus we can express $y_0^k$ as a 
linear combination of $y_i^a,\ i\geq 1, 1\leq a\leq N.$}
$y_0^k =- r/2 + \ldots,$ where the dots stand for a linear combination of 
$y_i^a,\ i\geq 1, 1\leq a\leq N,$ it expands (for each integer $r\in \Z$) 
as a power series in $\y$ (with $y_{0}^{k}$ excluded) with coefficients 
which are Laurent series in $\gl^{-1}$ 
(whose coefficients are {\em differential operators in $x$}
depending on $\x$ via $\T$, its translations and partial derivatives). 
The regularity condition means that all coefficients in front of the negative 
powers of $\gl$ vanish, i.e., the Laurent series are polynomials in $\gl.$

\begin{theorem}\label{t2} Let $\tau\in M$ be a semi-simple point. Then 
the total ancestor potential $\A_\tau^M$  satisfies the
HQE \eqref{HQE:ancestors}
\end{theorem}

\medskip

The HQE \eqref{HQE:ancestors} admit some kind of a classical limit. More 
precisely, the map $\tau_i\mapsto T_i$, $1\leq i\leq N,$ $Qe^{\tau_N}\mapsto T_{N+1}$,
identifies the rings $\C[\tau_1,\ldots,\tau_{N},Qe^{\tau_N}]$ and
$\C[T_1,\ldots,T_{N+1}].$ Given an element 
$f\in \C[\tau_1,\ldots,\tau_{N-1},Qe^{\tau_N}]$, we define the classical limit
of $f$ to be $f(0)$ where we are identifying $f$ with a polynomial in 
$\C[T_1,\ldots,T_{N+1}]$ and then we are setting $T_1=\ldots=T_{N+1}=0.$ Slightly
abusing the notations we will also say that we are setting 
$\tau_1=\ldots=\tau_N=Qe^{\tau_N}=0.$ The coefficients of the vertex
operators $\Gamma_\tau^a$ and $\Gamma_\tau^\gd$ depend polynomially on 
$\tau_1,\ldots,\tau_{N},Qe^{\tau_N}.$ After taking the classical limit
we obtain another set of Hirota quadratic equations. For more details
, see Sections \ref{sec:classical limit} and 
\ref{sec:From ancestors to descendents}. Here we summarize the answer. Put
\beq\label{vector:phi}
\f_\infty^\phi = 
\sum_{n\geq 0} \frac{\gl^n}{n!}d\tau_k(-z)^{-n-1}. 
\eeq
Denote by ${\rm Log}\, \gl$ a branch of the logarithmic function near 
$\gl=\infty.$ 
For each $1\leq a\leq k$, we introduce a vector in $\H:$ 
\beq\label{vector:a}
\f_\infty^a  = 
\frac{1}{k}{\bf g}_\infty^a 
+\sum_{i=1}^{k-1}
\sum_{n\in\Z}
\frac{\prod_{l=-\infty}^{n}(i/k - l) }
{\prod_{l=-\infty}^{0}(i/k - l)  }\, \gl^{i/k -n-1}\d_i (-z)^n ,
\eeq
where 
\ben
{\bf g}_\infty^a = 
\sum_{n\geq 0} \frac{\gl^n}{n!}(\log \gl - C_n)d\tau_k(-z)^{-n-1}+
\sum_{n\geq 0} n!\gl^{-n-1} d\tau_k z^n,
\een 
and $C_0:=0$, $C_n:=1+1/2+\ldots+1/n$ are the harmonic numbers. In 
the formulas above $a$ parametrizes different choices of $k$-th root of 1:
\ben
\log \gl = {\rm Log}\,\gl + 2\pi i (a-1),\quad  
\gl^{1/k} = \exp \(\frac{1}{k}\log \gl\).
\een
We also introduce a vector in $\H$ for 
each $b$ with $k+1\leq b\leq k+m:$ 
\beq\label{vector:b}
\f_\infty^b  =
-\frac{1}{m}{\bf g}_\infty^b 
-\sum_{j=1}^{m}
\sum_{n\in\Z}
\frac{\prod_{l=-\infty}^{n}(j/m - l) }
{\prod_{l=-\infty}^{0}(j/m - l)  }\, \gl^{j/m -n-1}\d_{k+m-j} (-z)^n ,
\eeq
where
\ben
{\bf g}_\infty^b = 
\sum_{n\geq 0} \frac{\gl^n}{n!}\[\log (\gl Q^{-m} ) - C_n\]
d\tau_k(-z)^{-n-1}+
\sum_{n\geq 0}n!\gl^{-n-1} d\tau_k z^n. 
\een
Just like above, $b$ parametrizes different choices of $m$-th root of 1:
\ben
\log \gl = {\rm Log}\,\gl + 2\pi i (b-k-1),\quad  
\gl^{1/m} = \exp \(\frac{1}{m}\log \gl\).
\een
Furthermore, introduce a vertex operator (acting on $\B_H$):
\ben
\Gamma_\infty^\delta = 
\exp\left( 
(\f_\infty^{\phi}-w_\infty)\ge\d_x \right)\sphat\ 
\exp \left({x}v_\infty/\ge \right)\sphat. 
\een
Finally, put
\ben
c_\infty^a = \frac{1}{k}\gl^{(1-k)/k},\ 1\leq a\leq k, \quad 
c_\infty^b = -\frac{1}{m}Q\gl^{-(1+m)/m},\ k+1\leq b\leq k+m.
\een
The limit of \eqref{HQE:ancestors} has the following form. 
\begin{definition}
{\em
We say that $\T\in B_H$ satisfies the HQE 
\eqref{HQE:descendents} if the 1-form
\beq
\label{HQE:descendents}
\(\Gamma_\infty^{\gd\#}\tensor \Gamma_\infty^\gd \) 
\(
\sum_a c_\infty^a\Gamma_\infty^{a}\tensor \Gamma_\infty^{-a} +
\sum_b c_\infty^b\Gamma_\infty^{b}\tensor \Gamma_\infty^{-b}  
\)
\(\T\tensor \T\)\  {d\gl}
\eeq
computed at $\q'$ and $\q''$ such that
$\widehat w_\infty' - \widehat {w}_\infty'' =r$, is regular in 
$\gl$ for each $r\in \Z$.
}
\end{definition}
We remark that \eqref{monodromy:killer} holds with $\tau=\infty,$ which implies that
the expression \eqref{HQE:descendents} is single-valued near $\gl=\infty$ 
and independent of our choice of the branch ${\rm Log}\,\gl.$ The regularity 
condition is interpreted as before. 
 
\begin{theorem}\label{t1} 
The total descendent potential $\D^M$ satisfies the
HQE \eqref{HQE:descendents}.
\end{theorem}

\subsection*{Acknowledgments}
We thank Y. Ruan for his interests in this work. Part of this work was pursued 
in the Mathematical Sciences Research Institute where the second author held a 
postdoctoral fellowship in the spring of 2006. It is a pleasure to acknowledge 
its hospitality and support. Finally, we want to thank the referee for pointing
out some inaccuracies and fixing a numerous number of misprints.


\section{Hirota Quadratic Equations and $A_1$ singularity}\label{a1kdv}

According to the main result in \cite{G2}, Theorems \ref{t2} and \ref{t1}
have analogues for the Frobenius manifold $M_{A_1}$ corresponding to $A_1$ 
singularity. Here we review the results. $M_{A_1}\simeq \com$ is the space of quadratic polynomials
$f = x^2/2 + u,\ u\in \C.$ The metric and the Frobenius multiplication are 
the standard ones of $\C.$ The construction of a symplectic loop space $\H$ and 
corresponding Bosonic Fock space $B$ depends only on a vector space equipped with 
a non-degenerate bi-linear form. Let $\H_\C$ and $B_\C$ be the corresponding 
objects for $\C$ equipped with the standard pairing. 
 
Let $x_\pm = \pm\sqrt{2(\gl-u)}$ be the 
solutions of $f(x)=\gl.$  Then we define period vectors
\ben
I_\pm^{(-n)}(u,\gl) & := & -\d_u \frac{1}{n!} \int_{x_\pm} d^{-1}\(\(\gl-f\)^n dx\)= 
\pm\frac{[2(\gl-u)]^{n-1/2}}{(2n-1)!!}, \\
I_\pm^{(n)}(u,\gl) & := & \d_\gl^n I_\pm^{(0)}(u,\gl) =
\pm(-1)^n \frac{(2n-1)!!}{[2(\gl-u)]^{n+1/2}} . 
\een
Let $\Gamma_u^\pm$ be the vertex operator (acting on $B_\C$) corresponding to
the vector
\ben
\f_u^\pm = \sum_{n\in \Z} I_\pm^{(n)}(u,\gl)(-z)^n. 
\een
It is known that the potentials $\D^{A_1}$ and $\A_u^{A_1}$ coincide with
the so called {\em Witten--Kontsevich tau-function}: 
\ben
\D_{\rm pt}(\t) = \exp \left(\sum_{g,n} \frac{\ge^{2g-2}}{n!} 
\int_{\overline{\M}_{g,n}} \t(\psi_1)\ldots \t(\psi_n) \right),
\een  
where 
$\t=(t_0,t_1,\ldots)$ is a sequence of formal variables, 
$\overline{\M}_{g,n}$ is the moduli space of stable Riemann surfaces of
genus $g$ with $n$ marked points, 
$\psi_i$ is the first Chern class of the $i$-th universal cotangent line bundle on $\overline{\M}_{g,n}$,
$\t(\psi_i) = \sum_l t_l \psi_i^l,$  and the sum is over all $g$ and $n$ 
with the convention that the integral is $0$ if  $\overline{\M}_{g,n}$
is empty. 
$\D_{\rm pt}$ is identified with an element of the Fock space $B_\C$ via the 
{\em dilaton shift} $\t(z) = \q(z)+z.$ 

According to \cite{G2}, corollary of Proposition 2, 
the Witten--Kontsevich tau-function satisfies the following HQE: the 1-form
\beq \label{HQE:kdv}
\(
\frac{1}{\sqrt{2(\gl-u)}}\, \Gamma_{u}^{+} \tensor \Gamma_{u}^{-} - 
\frac{1}{\sqrt{2(\gl-u)}}\, \Gamma_{u}^{-} \tensor \Gamma_{u}^{+} \)
\(\D_{\rm pt} \tensor \D_{\rm pt} \)d\gl.
\eeq 
is regular in $\gl$ (in the sense explained in the Introduction). 

We make several remarks. First of all the coefficients in front
of the vertex operators are precisely $c_f^\pm(\gl) = 1/f'_x(x_\pm)$
which agrees with the formula for $c_f^i(\gl)$ in the case $f\in M$ -- the
space of Laurent polynomials. Second, when $u=0,$ \eqref{HQE:kdv} is 
precisely the Witten's conjecture \cite{W}, proved by Kontsevich 
\cite{Ko}. Finally, the proof of \eqref{HQE:kdv} for $u\neq 0$ 
follows from the case $u=0$ and {\em the string equation}.


\section{Frobenius structure on the space of Laurent polynomials}
\label{froblaurent}


\subsection{Flat structure}  
\label{sec:flat_structure}

In this section, following \cite{DZ}, we will show that the residue
metric \eqref{res_metric} is flat. 

Let $t=(t_1,t_2,\ldots,t_N)$ and $f_t$ be the corresponding Laurent polynomial in $M.$ 
In a neighborhood of $\gl=\infty$ the equation 
$f_t(x)=\lambda$ has two types of solutions depending on whether
$x$ is close to $\infty$ or to $0$. Let $x_a,\ a=1,2,\ldots, k$ be 
the solutions close to $\infty$ and $x_b,\ b=k+1,k+2,\ldots, k+m$ the 
solutions close to $0.$ They expand as series in $\gl$ as follows:
\ben
&&
\log x_a = \frac{1}{k}\[\log \gl - \tau_1 \gl^{-1/k} - \ldots - 
\tau_{k-1}\lambda^{-(k-1)/k} - \tau_k\gl^{-1}\] + O\(\gl^{-1-1/k}\), \\
&&
\log x_b = \frac{1}{m}\[-\log \frac{\lambda}{Q^m} + \tau_{k+m} + 
\tau_{k+m-1} \gl^{-1/m} + \ldots + 
\tau_k \gl^{-1}\] + O\(\gl^{-1-1/m}\),
\een
where indices $a$ and $b$ have the same ranges as above. They parametrize 
different choices of $k$-th and $m$-th root of $1$ respectively. 
The coefficients $\tau_i$, $i=1,2,\ldots, k+m$ can be expressed in terms 
of $t_1,\ldots, t_N$ as follows:
\ben
&&
\tau_i = -\frac{k}{i} \res_{x=\infty} \ f_t(x)^{i/k}\,\omega ,\quad 
1\leq i\leq k-1  \\
&&
\tau_{k+m-j}= \frac{m}{j} \res_{x=0} \ f_t(x)^{j/m}\, \omega ,\quad
1\leq j\leq m-1 \\
&& \tau_N = mt_N, \quad \tau_k = t_k.
\een
Using these formulas we get
\beqa\label{flat_coord}
&&
t_i = \tau_i + f_i(\tau_1,\ldots, \tau_{i-1}),\quad 1\leq i\leq k-1\\
&& \label{coord_flat}
t_{k+j}= \tau_{k+j} + 
h_j(\tau_{k+j+1},\ldots,\tau_{k+m-1}),\quad 1\leq j\leq m-1 
\eeqa
where $h_j$ and $f_i$ are certain polynomials of degrees at least 2. Thus 
the corresponding Jacobian is non-degenerate and the functions $\tau_1,\ldots,\tau_{N-1},e^{\tau_N}$
give a coordinate system on $M.$  
Moreover, according to \cite{DZ}, in such coordinates the residue metric has the form 
\ben
&& (\d/\d{\tau_i}, \d/\d{\tau_{k-i}})_\tau = 1/k,\quad
i=1,2,\ldots, k-1 \\
&&
(\d/\d{\tau_{k+j}}, \d/\d{\tau_{k-j+m}})_\tau =1/ m,\quad
j=0,1,\ldots, m, 
\een 
and all other pairings between $\d/\d\tau_i,1\leq i\leq N$ are 0.  

For the sake of completeness let us show how
to compute $(\d/\d{\tau_i}, \d/\d{\tau_{k-i}})_\tau.$ Put $\xi = \gl^{1/k}$
and identify $\xi$ with a new coordinate on $\C^*$ near $x=\infty,$ related 
to $x$ via $f_\tau(x)=\xi^k.$ By chain rule we have 
$
\d_{ \tau_i}f_\tau + \d_x f_\tau\d_{\tau_i} x = 0,
$
which implies that $\d_{\tau_i} f_\tau \omega = -(\d_{\tau_i} \log x)df_\tau.$ 
Using the expansion of $\log x$ from above, we find  
$
(\d/\d\tau_i)(f_\tau\omega) = k\xi^k [k^{-1}\xi^{-i} + O(\xi^{-k}) ]{d\xi}/{\xi}.   
$
Note that in the residue pairing only the residue at
$x=\infty$ contributes:
\ben
(\d_{\tau_i},\d_{\tau_j})_\tau = 
- \res_{\xi=\infty}\  [\frac{1}{k}\xi^{k-i-j} + O(\xi^{-1})]d\xi/\xi=
 \frac{1}{k}\,\delta_{i+j,k}.
\een 
In flat coordinates the Euler vector field takes on the form:
\ben
E= \tau_k\d_{\tau_k} + 
\sum_{i=1}^{k-1} \frac{i}{k}\tau_i\d_{\tau_i} + 
\sum_{j=1}^{m-1} \(1-\frac{j}{m}\)\tau_{k+j}\d_{\tau_{k+j}}+ 
\(\frac{1}{k}+\frac{1}{m}\)m\d_{\tau_{N}}.
\een

\subsection{Oscillating integrals}

Let $\{\tau_i\}_{i=1}^N$ be the flat coordinates introduced above.
Denote the corresponding coordinate vector fields by $\d_i:=\d/\d \tau_i.$ 
The following lemma is crucial for our construction. Probably it could be
derived from \cite{DS1,DS2} or \cite{Ba}. However we prefer to give a direct proof.
\begin{lemma}\label{lemma:primitive_form}
For each $i$ and $j$, $1\leq i,j\leq N$ there is a Laurent polynomial 
$G_{ij}(\tau,x)$ in $x$ such that
\beq\label{primitive form}
\frac{\d^{2}f_\tau}{\d\tau_i\d\tau_j}\omega = dG_{ij}, \quad
\((\d_if_\tau)\, (\d_jf_\tau) - \sum_p A_{i,j}^p \d_pf_\tau\)\omega = 
G_{ij}\, df_\tau,  
\eeq
where $d$ is the De Rham differential on $\C^*.$ 
\end{lemma}
\proof
From the definition of $\bullet$, it follows that the second equality holds 
for a uniquely determined Laurent polynomial 
\ben
G_{ij} = \sum_{a=-m}^{k-1} G_{ij}^{(a)}(\tau) x^a.
\een
Write the polynomial $f_\tau$ as 
$x^k+T_1 x^{k-1} + \ldots + T_N x^{-m}.$ Then we need to show that 
\ben
\frac{\d^{2}T_a}{\d\tau_i\d\tau_j} = (k-a)G_{ij}^{(k-a)},\quad 1\leq a\leq N.
\een
Assume first that $1\leq a\leq k.$ Let $\xi$ be a new coordinate in a 
neighborhood of $x=\infty$ defined by $f_\tau(x) = \xi^k,$ and so $\log x$ 
is expressed in terms of $\xi$ according to the expansions defining the 
flat coordinates $\tau_i,\ 1\leq i\leq k,$ except that we need to put $\gl=\xi^k.$  
Then we have
\ben
T_a & = & 
-\res_{x=\infty} f_\tau(x) x^{a-k-1} dx = 
-\frac{1}{a-k}\res_{\xi=\infty} \xi^k d\( x^{a-k} \) \\
&=&
 \frac{k}{a-k} \res_{\xi=\infty} x^{a-k} \xi^{k-1}d\xi,
\een
where $x=\xi(1 + O(\xi^{-1})).$ Using this formula we compute $\d_i\d_j T_a:$ 
\ben
\frac{\d^{2}T_a}{\d\tau_i\d\tau_j} &=& 
k \res_{\xi=\infty}\[ (a-k-1) x^{a-k-2}(\d_ix)(\d_jx) + 
x^{a-k-1}\d_i\d_jx\] \xi^{k-1}d\xi \\
&=& 
k \res_{\xi=\infty}\[ (a-k) x^{a-k} (\d_i\log x)(\d_j\log x) + 
                            x^{a-k}  \d_i\d_j\log x \]
\xi^{k-1}d\xi.
\een
Notice that the last term in the square brackets does not contribute to the residue
because the highest possible power of $\xi$ is $a-k -k-1$. Therefore, after
passing back to the old coordinate $x$ and using that $(\d_i\log x) df=-\d_if\omega$, we get
\ben
\frac{\d^{2}T_a}{\d\tau_i\d\tau_j} & = &
(a-k)\res_{x=\infty}\[x^{a-k-1}\, {(\d_i f_\tau)(\d_jf_\tau)}\frac{\omega}{f_\tau'}\]  \\
& =&
(a-k)\res_{x=\infty}\[x^{a-k}G_{ij}+  
x^{a-k-1}\(\sum_p A_{ij}^p \d_p f_\tau\)\frac{1}{f_\tau'} \] \omega.
\een
Again, the second term in the square brackets does not contribute to the residue
because the highest possible power of $x$ is ${a-k-1} \leq -1.$
The residue of the first term is clearly $-G_{ij}^{(k-a)}.$ 

In the case when $k+1\leq a\leq N$, we  
pass to a new coordinate $\eta$ near $x=0$ via $f_\tau(x) = T_N \eta^{-m}.$ 
Then proceed by a similar argument. 
\qed

In particular,  the oscillating integrals 
\beq\label{osc}
\J_{\B}(\tau,z) = (-2\pi z)^{-1/2}\int_\B e^{f_\tau/z}\omega
\eeq
satisfy the following differential equations
\beq
\label{frob_str}
z\frac{\d^2 \J_\B}{\d \tau_i\d \tau_j} = 
\sum_{p=1}^{N} A_{ij}^p\frac{\d\J_\B}{\d \tau_p},
\eeq
where the integration cycle $\B$ is an element of the relative 
homology group 
\beq\label{homology_osc}
\lim_{M\rightarrow \infty} 
H_1(\C^*, \{x\in\C^*\  |\ {\rm Re }({f_\tau/z})<-M\};\Z)\iso \Z^N.
\eeq
The oscillating integral $\J_\B$ also satisfies some homogeneity conditions
due to the fact that $f_\tau$ and $\omega$ are homogeneous:
\beq\label{frob_e}
\(z\d_z + E+1/2\) \J_\B = 0. 
\eeq
Let $J_\B$ be a vector field on $M$ defined by 
\ben
(J_\B(\tau,z),\d_i)_\tau = z\d_i\,\J_\B.
\een
Then equations \eqref{frob_str} and \eqref{frob_e} are equivalent to 
\beq\label{de_J}
z\d_{i} J_\B = (\d_{i}\bullet)\, J_\B,\quad 
\(z\d_z + E\) J_\B = \mu J_\B,
\eeq
where $\mu$ is {\em the Hodge grading operator}: 
\ben
\mu(\d_i) = \(\frac{i}{k}-\frac{1}{2}\)  \d_i,\quad 
\mu(\d_{k+j})=\(\frac{1}{2}-\frac{j}{m}\) \d_{k+j}, \quad 
1\leq i\leq k-1,0\leq j\leq m.
\een  
Let $\B_i$ be a basis of cycles in the relative homology group 
\eqref{homology_osc}. Then the matrix $J$ with columns $J_{\B_i}$ 
is a fundamental solution to the system \eqref{de_J}. This means 
that we can extend the connection $\nabla$ defined in \eqref{flat_connection}
to a connection on the trivial bundle on $M\times\C^*$ with fiber $H$ by
setting 
\ben
\nabla_{\d/\d z} := \d/\d z + \(E\bullet \)z^{-2} - \mu z^{-1}.
\een
The extended connection is flat because the corresponding system of
differential equations admits a fundamental solution. 

For each $\tau\in M$, the $z$-direction of $\nabla$ defines a 
connection on $\C \proj^1$ which has an irregular singular point at
$z=0$ and a regular singular point at $z=\infty.$  At $z=0$ the fundamental 
solution of $\nabla_{\d/\d z} \Phi=0$ admits a certain asymptotic and 
at $z=\infty,$ $\nabla_{\d/\d z}$ can be transformed
via a gauge transformation into a {\em canonical form.} These two
ingredients, the asymptotic and the gauge transformation, contain the
essential information about the Frobenius structure. They will
be used to define the total descendent and the total ancestor potentials 
of $M.$ 


\subsection{Stationary phase asymptotic}

Let $\tau\in M$ be {\em a semi-simple}  point, i.e., $f_\tau$ has  only 
Morse type critical points $q_i,\ 1\leq i \leq N.$ Denote the 
corresponding critical values by $u_i.$ They form a coordinate system
called {\em canonical coordinate system}. Let $\Delta_i$ be the 
Hessians of $f_\tau$ at $x_i$ with respect to the volume form $\omega.$ Then 
the linear map 
\ben
\Psi: \C^N\rightarrow T_\tau M,\quad 
\Psi(e_i):={\bf 1}_i:=\sqrt{\Delta_i}\d/\d u_i, 
\een
is an isomorphism of Frobenius algebras. Here $e_i,\ 1\leq i\leq N$ are 
the standard coordinate vectors and a Frobenius algebra structure on  $\C^N$ 
is defined by the product: $e_ie_j = \delta_{ij}\sqrt{\Delta_i}\,e_i$ and 
the metric:  $(e_i,e_j)=\delta_{ij}$. 

Furthermore, choose a basis of cycles $\B_i$ in \eqref{homology_osc} (e.g., 
by means of Morse theory for ${\rm Re}(f_\tau/z)$) and   
let $J:\C^N\rightarrow T_\tau M$ be a linear operator defined by 
$J(e_i) = J_{\B_i}.$ Then there is an asymptotical expansion $J\sim \Psi R e^{U/z}$ as 
$z\rightarrow 0,$ valid in some sector in the $z$-plane, such that $U$ is a diagonal 
matrix with entries $u_1,\ldots, u_N$ and $R=1+R_1z+\ldots$ is a certain
series with matrix coefficients.  Such an asymptotic will be derived in the proof of \leref{vector}. Thus the system 
\eqref{de_J} admits an asymptotical solution of the form $\Psi R e^{U/z},$
which {\em a priori} might depend on the choice of the cycles $\B_i$.  
However, according to \cite{G1}, Proposition, part (d), such a solution
is unique and it automatically satisfies the symplectic condition 
$R^*(-z)R(z)=1.$


\subsection{Calibration of $M$}

\begin{proposition} There exist a gauge transformation  
$S_\tau=1+S_1z^{-1}+\ldots,$ $S_p\in {\rm End}(H)$ satisfying the 
symplectic condition $S_\tau^*(-z)S_\tau(z)=1,$ such that 
\beq\label{de_S}
z\d_{i} S_\tau = (\d_i\bullet_\tau)\,S_\tau ,\quad 
z\d_z S_\tau = [\mu, S_\tau]-\((E\bullet)S_\tau - S_\tau\rho\)z^{-1}, 
\eeq
where $\mu$ is the Hodge grading operator and $\rho$ is the cup product 
multiplication by $(1/k+1/m) m\d_N .$
\end{proposition}
\proof
The first equation in \eqref{de_S} gives us the following recursive 
relation: $\d_i S_p = (\d_i\bullet)S_{p-1}.$ The multiplication 
operator $\d_i\bullet$ depends polynomially on 
$\tau_1,\ldots,\tau_N,$ and  $Qe^{\tau_N}.$ Thus starting from $S_0=1$ we 
can recover uniquely all other $S_p,\ p\geq 1$ by integrating the 
recursive relations and requiring that $S_p$ vanishes when 
$\tau_1,\ldots,\tau_N, Qe^{\tau_N}$ are set to $0.$

We claim that the so constructed series $S$ automatically satisfies 
the second equation of \eqref{de_S} and the symplectic condition. 
Identify $S(\tau,z)$ with a section of the bundle $\pi^*\({\rm End} TM\),$
where $\pi:M\times\C^*\rightarrow M$ is the projection.  
Then $\eqref{de_S}$ means that $S$ is a horizontal section of the 
following connection
\ben
\nabla^{\rm End} = d- \sum_{i=1}^N (\d_i\bullet)z^{-1}d\tau_i - 
 \({\rm ad}(\mu)z^{-1} - (E\bullet -\rho^R)z^{-2} \)dz,
\een
where $\rho^R$ is the classical multiplication by $(1/k+1/m)m\d_N$ 
{\em from the right.}
The flatness of $\nabla$ implies the flatness of $\nabla^{\rm End}$. 
In particular, if we set $\Phi(\tau,z) = \nabla^{\rm End}_{\d/\d z} S$
then $\nabla^{\rm End}_{\d/\d \tau_i }\Phi = 0.$ However, $\Phi$ is a
power series in $z^{-1}$ with coefficients depending polynomially
on  $\tau_1,\ldots,\tau_N, Qe^{\tau_N}$, and $\Phi$ vanishes when those variables 
are set to $0$. Thus $\Phi=0$, which is precisely
the second equation in \eqref{de_S}. 

Let us  prove that $S$ satisfies the symplectic condition 
$S^*(-z)S(z) =1.$ Differentiate with $\d_i$ and use the first equation
in \eqref{de_S} and the fact that the operators $\d_i\bullet$ are self 
adjoint, we get that $\d_i \( S^*(-z)S(z)\)=0,$ i.e., $S^*(-z)S(z)$
is a constant independent of $\tau.$ Set 
$\tau_1,\ldots,\tau_N, Qe^{\tau_N}$ to $0$, then $S=1$ by construction. Thus $S^*(-z)S(z) =1.$  
\qed

The choice of a solution to \eqref{de_S} is called {\em a calibration}.
We choose a calibration of $M$ as follows: each coefficient $S_p$ depends
polynomially on $\tau_1,\ldots,\tau_N,$ and $Qe^{\tau_N}.$ We require that
$S_\tau = 1$ when  the variables $\tau_1,\ldots,\tau_N, Qe^{\tau_N}$ 
are set to 0.


\subsection{Descendents and ancestors}
\label{sec:desc_anc}

By definition, the {\em twisted loop group} is  
$$
\L^{(2)}\lieGL(H) = \left\{
M(z)\in \L\lieGL(H)\ |\ M^*(-z)M(z)=1\right\},
$$
where $*$ means the transposition with respect to the bilinear pairing. 
The elements of the twisted loop group of the form
$M=1+M_1z+M_2z^2+\ldots$ (respectively $M=1+M_1z^{-1}+M_2z^{-2}+\ldots)$
are called upper-triangular (respectively lower-triangular) linear 
transformations. 
They can be quantized as follows: write $A=\log M,$ then $A(z)$ is an
infinitesimal symplectic transformation. We define  
$\widehat M=\exp \hat A,$ where $A$ is identified with the 
quadratic Hamiltonian $\Omega(A\f,\f)/2$ and on the space of
quadratic Hamiltonians the quantization rule $\sphat\ $ is defined by:
\ben
(q_{k,i}q_{l,j})\sphat := \frac{q_{k,i}q_{l,j}}{\ge^2},\ \ 
(q_{k,i}p_{l,j})\sphat := q_{k,i}\frac{\d}{\d q_{l,j}},\ \ 
(p_{k,i}p_{l,j})\sphat := \ge^2\frac{\d^2}{\d q_{k,i}\d q_{l,j} }.
\een 
We remark that $\sphat\ $ defines only {\em a projective representation} of
the lower- and upper-triangular subgroups.

Motivated by Gromov--Witten theory, A. Givental \cite{G1} introduced the 
so-called {\em total ancestor} and {\em total descendent} potentials of 
a semi-simple Frobenius manifold. 
The total ancestor potential $\A_\tau$ of $M$ is defined for any semi-simple
point $\tau\in M$ as follows:
\ben
\A_\tau = \widehat\Psi\, \widehat R\,  \(e^{U/z}\)\sphat\, \prod_{i=1}^N 
\D_{\rm pt}^{(i)}\Delta_i^{-1/48},
\een   
where $\D_{\rm pt}$ is the Witten--Kontsevich tau-function, $\D_{\rm pt}^{(i)}$
is a vector in the Fock space $B_{\C^N}$ defined by 
$\D_{\rm pt}^{(i)}(\q)=\D_{\rm pt}(\q^i)$ with $\q=\sum \q^i e_i \in \C^N[z].$ 
The linear operators $R$ and $e^{U/z}$
are elements of the upper- and lower-triangular twisted 
loop subgroups respectively. They act on $B_{\C^N}$ according to the above quantization rules. 
Finally, $\widehat\Psi$ is just the identification between the 
Fock spaces $B_H$ and $B_{\C^N},$ i.e., $(\widehat \Psi \mathcal{G})(\q) = 
\mathcal{G}(\Psi^{-1}\q).$
   
The total descendent potential is defined by
\ben
\D = {C(\tau)}\widehat S_\tau^{-1} \A_\tau, \quad 
C(\tau)=\exp \(\frac{1}{2}\int^\tau R_1^{ii}du_i\),
\een
where $R_1^{ii}$ are the diagonal entries of $R_1\in\lieGL(\C^N).$ 
The constant $C$ is chosen in such a way that $\D$ is independent of $\tau.$

\section{$\mathbb{P}^1$-orbifolds}\label{p1orb}
\label{sec:ocp1}


In this Section we discuss Gromov-Witten theory of the orbifold $\mC_{k,m}$-- 
an orbifold obtained from $\mathbb{P}^1$ by cutting two small discs
$D_1=\{|z|\leq \ge\}$ and  $D_2=\{|z^{-1}|\leq \ge\}$ respectively near 
$z=0$ and $z=\infty$ and gluing back the orbifolds $D_1/\Z_k$ and $D_2/\Z_m$ in the
obvious way. The main goal is to compare the Frobenius manifold $M_{k,m}$ and a Frobenius manifold corresponding to the orbifold quantum cohomology of $\CC_{k,m}$. Our approach is to compute the {\em small} orbifold quantum cohomology of $\CC_{k,m}$ and then use a reconstruction result, see Theorem \ref{reconstruction}.
 
 \subsection{Reconstruction theorem}

Let $M$ be a small ball centered at $0$ in $\C^N.$ Assume that $g$ is a non-degenerate bi-linear pairing on $TM$, $A$ is a holomorphic section of $T^*M^{\tensor 2}\tensor TM,$ i.e., the tangent spaces $T_tM$ are equipped with a multiplication $\bullet_t$ which depends holomorphically on $t\in M,$  $e$ is a vector field on $M$ such that its restriction to $T_tM$ is a unity with respect to $\bullet_t,$ and finally $E$ is a vector field on $M.$ 

\begin{definition}\label{def:frob}
{\em
The data $(M,g,A,e,E)$ form a Frobenius structure on $M$ if the following conditions are satisfied. 
\begin{enumerate}
\item
$g$ and $\bullet$ satisfy the Frobenius property: 
$g(X\bullet Y_1,Y_2) = g(Y_1,X\bullet Y_2),$

\item
The one-parameter group corresponding to $E$ acts on $M$ by conformal
transformations of $g$, i.e., $\L_Eg =Dg,$ for some constant $D\in \C,$ 

\item
$e$ is a flat vector field: $\nabla^{\rm L.C.} e=0,$ where $\nabla^{\rm L.C.}$
is the Levi-Civit\'a connection of $g,$

\item
The connection operator    

\beq\label{str_connection}
\nabla = \nabla^{\rm L.C.} - {z^{-1}} \sum_{i=1}^N \(\frac{\d}{\d {t_i}}\bullet_t\)dt_i 
+ \(z^{-2} (E\bullet_t) - z^{-1}\mu \)dz,  
\eeq
where $\mu := \nabla^{\rm L.C.} (E)-(D/2){\rm Id}:TM\rightarrow TM$ 
is {\em the Hodge grading operator,} is flat, i.e., $\nabla^2=0.$
\end{enumerate}
}
\end{definition}
Here $\{t_i\}$ are arbitrary coordinates on $M$, $\d_{t_i}\bullet_t$ (respectively $E\bullet_t$) is the $\bullet_t$-multiplication by the vector field $\d_{t_i}$ (respectively $E$), and $\nabla$ is a connection on the bundle $\pi^*(TM)$ with base $M\times \C^*,$ where $\pi:M\times \C^*\rightarrow M$ is the projection.  

Let us assume that $0\in M$ is a semi-simple point, i.e., the Frobenius algebra $T_0M$ is diagonalizable. Equivalently, there are local coordinates $u^i,1\leq i\leq N$, called {\em canonical coordinates} which diagonalize the metric $g$ and the multiplication $\bullet_\tau$:
\ben
\d/\d u^i\bullet \d/\d u^j= \delta^{ij}\d/\d u^j,\quad 
g(\d/\d u^i, \d/\d u^j)=\delta^{ij} \theta_j,\quad 1\leq i,j\leq N,
\een
where $\theta_j,$ are some holomorphic functions on $M$. Moreover, from the flatness of the connection operator \eqref{str_connection}, it follows that the coordinates $u^i$ could be chosen such that the Euler vector field assumes the form $E=\sum_iu^i\d/\d u^i$ (see \cite{D}, Lemma 3.5). The goal in this subsection is to prove the following theorem:

\begin{theorem}\label{reconstruction}
Let $M$ be a holomorphic Frobenius manifold and $0\in M$ is a semi-simple point such that the following conditions are satisfied:
\begin{enumerate}
\item\label{kth_root} 
The restriction $E_0$ of the Euler vector field to $0$ has a $k$-th root $v\in T_0M$ (i.e. $E_0=v^k$ ) such that $v$ is invertible and $v$ generates the Frobenius algebra $T_0M.$ 
\item\label{spectral_condition}
Let $\overline{\mu}$  and ${\rm Diag}\(u^1_{(0)},\ldots,u^N_{(0)}\)$
be the matrices respectively of the Hodge grading operator $\mu$ and 
the operator of multiplication by the Euler vector $E_0$, 
in a basis of $T_0M$ which diagonalizes $\bullet_0.$ If $(i,j)$ is a 
pair of indices such that $u^i_{(0)}=u^j_{(0)}$  then the $(i,j)$-th
entry of $\overline{\mu}$ is zero. 
\end{enumerate}
Then the Frobenius structure on $M$ is uniquely determined from the Frobenius
algebra $T_0M$ and the Hodge grading operator $\mu$. 
\end{theorem}

The idea of the proof is to reconstruct successively the terms  of the Taylor's  expansion of $u^i$, $1\leq i\leq N.$ We use certain recursive relations, constructed from condition \eqref{kth_root} and the flatness of \eqref{str_connection}. Condition \eqref{spectral_condition} guarantees that the recursive relations can be solved. Our argument was inspired by the proof of Lemma 2.9 in \cite{HM}. 

\proof
Condition (4) in \deref{def:frob} implies that $\nabla^{\rm L.C.}$ is a flat connection and that the multiplication $\bullet_t$ is associative and commutative. Let us identify $T_0M$ with $\C^N$ by fixing a basis $\phi_1:=e,\phi_2:=v,\phi_3,\ldots,\phi_N$ of $T_0M.$ Using the flat connection $\nabla^{\rm L.C.}$ we extend $\phi_a,1\leq a\leq N$ to vector fields $\d_a:=\d/\d \tau^a$ on $M,$ where $\tau=(\tau^1,\ldots,\tau^N)$ is a flat coordinate system on $M,$ and so all other tangent spaces $T_\tau M$ are canonically identified with $\C^N$ as well. 

Let $U:= E\bullet_\tau$ and $F_a:= \d_a\bullet_\tau$, $1\leq a\leq N$ be the linear operators of multiplication by the corresponding vector fields. In view of the above identifications, we may regard $U$ and $F_a$ as $N\times N$-matrices, whose entries are holomorphic functions in $\tau.$ We define a grading on the space of holomorphic $N\times N$- matrices by assigning degree 1 to each of the coordinate functions $\tau^1,\ldots,\tau^N.$  If $A(\tau)$ is a holomorphic matrix then we denote by $A^{(0)}+A^{(1)}+A^{(2)}+\ldots$ its homogeneous decomposition, i.e., $A^{(n)}$ is a finite sum of matrices whose entries are monomials of degree $n$. For some matrices, in order to avoid cumbersome notations, we write $A_{(n)}$ instead of $A^{(n)}$. Also, we denote by $A^{(\geq n)}$ (resp. $A^{(\leq n)}$) the matrix obtained from $A$ by truncating all terms of degree $<n$ (resp. $>n$).

Let us denote by $\Psi$ the matrix whose $(i,a)$-entry (i.e. $i$-th row and $a$-th column) is given by: $\Psi_{ia} = \d_a u^i.$  It is easy to see that: 
\beq\label{diagonalize_a}
\Psi F_a \Psi^{-1} = {\rm Diag} \Big(
\d_a u^1,\ldots,\d_a u^N \Big)=: D_a,
\quad 1\leq a\leq N
\eeq
and 
\ben
\Psi U \Psi^{-1} = {\rm Diag} (u^1,\ldots, u^N )=: D.
\een
We know $U^{(0)}$ and $F_a^{(0)}, 1\leq a\leq N$ and we want to reconstruct
$U^{(n)}$ and $F_a^{(n)}$ for $n>0.$ 
 
Note that the matrix $U$ admits a holomorphic $k$-th root $V.$ Indeed, we have $U=\Psi^{-1}\,D\,\Psi$ and $D$ is diagonal with entries 
$u^i(\tau) = u^i_{(0)}+u^i_{(\geq 1)}$. On the other hand $u^i_{(0)}$ is non-zero, because it is a $k$-th power of an eigenvalue of $v\bullet_{0}$ and the later is an invertible matrix by definition. Therefore we may define (by using the binomial formula):
\ben
(u^i)^{1/k}:= \(u^i_{(0)}\)^{1/k}
\Big( 1 + \sum_{j=1}^{\infty} \binom{1/k}{j} \(u^i_{(\geq 1)}/u^i_{(0)}\)^j
\Big).
\een  
Therefore $V:=U^{1/k}:= \Psi\,D^{1/k}\,\Psi^{-1}$ is a $k$-th root of $U.$ Moreover, without any restrictions we may assume that $V_{(0)}=v\bullet_{0}.$ Since $v$ generates the Frobenius algebra $T_0M$, we can find polynomials $f_a(x),\ 1\leq a\leq N$ such that $\phi_a = f_a(V_{(0)})e$, where $e\in H$ is the unity, or equivalently $F_a^{(0)}= f_a(V_{(0)}).$ 

Assume that we have determined the matrices $U^{(i)},V^{(i)},F_a^{(i)}$ for all $a=1,2,\ldots,N$ and all $i=0,1,\ldots,n-1$. We want to prove that the matrices for $i=n$ are uniquely determined as well. Let us remark that there is a small difference between the cases $n=1$ and $n>1$ which however appears only in the proof of \leref{lemma_B}, part b), below.

From the flatness of the connection operators \eqref{str_connection} we have $[\nabla_{\d_a},\nabla_{\d/\d z}]=0$. Comparing the terms of degree $n-1$ we get: 
\beq\label{eqn_f}
\d_a U^{(n)} = F_a^{(n-1)} +[\mu,F_a^{(n-1)}].
\eeq
A direct corollary of this equation is that $U^{(n)}$ is uniquely determined from $F_a^{(n-1)},1\leq a\leq N.$

Put $P:=\Psi^{(0)}.$ The entries of this matrix can be determined as follows. We know that $P$ diagonalizes the Frobenius product $\bullet_{0}$. Therefore we have:
\ben
P\,V_{(0)}\, P^{-1} = {\rm Diag}(\gl_1,\ldots,\gl_N)
\een
and 
\ben
P\, F_a^{(0)}\, P^{-1} =P\, f_a(V_{(0)})\, P^{-1}= 
{\rm Diag}(f_a(\gl_1),\ldots,f_a(\gl_N)),\quad
1\leq a\leq N.
\een
Comparing with \eqref{diagonalize_a} we get that the $(i,a)$-entry of $P$ is given by:
\ben
P_{ia} = f_a(\gl_i)=:\gl_{a,i}.
\een
Let us remark that the eigenvalues $\gl_{a,i}$ have the following two properties. First, according to our choice of a basis of $T_0M$ we have $\phi_1=e$ and $\phi_2=v$, therefore $\gl_{1,i}=1$ and $\gl_{2,i}=\gl_i$. Second, the eigenvalues $\gl_i,1\leq i\leq N$ are pairwise different. Indeed, if this is not the case then, there exists a non-diagonal matrix $A\neq 0$ that commutes with $V_{(0)},$ and hence it commutes with $F_a^{(0)}=f_a(V_{(0)}),$ $1\leq a\leq N$. This is impossible because, $A$, $F_a$,  $1\leq a\leq N$ are linearly independent and a maximal abelian Lie subalgebra of $\liegl(N,\C)$ has dimension $N.$   

Given a $N\times N$ matrix $A$ we put $\overline{A}=P\,A\, P^{-1}.$ 
Let us compare the degree $n$ terms in the equation $U=V^k$. We get
\beq\label{nexc}
\overline{U}^{(n)} = \sum_{s=0}^{k-1}\overline{V}_{(0)}^{s} \overline{V}_{(n)}
\overline{V}_{(0)}^{k-1-s} + \ldots, 
\eeq
where the dots stand for terms which depend on $V_{(i)}$ with $i\leq n-1.$
Note that the $(i,j)$ entry of the matrix sum from above is 
\beq\label{nexc2}
\[\overline{V}_{(n)}\]_{ij}\(\gl_i^{k-1}+\gl_i^{k-2}\gl_j+\ldots+\gl_j^{k-1}\).
\eeq
The above sum of $\gl$'s is zero precisely when the pair $(i,j)$ is 
such that $\gl_i\neq\gl_j$ but $\gl_i^k=\gl_j^k.$ We call such a pair 
{\em exceptional} and the entries in a $N\times N$ matrix corresponding
to an exceptional pair are called exceptional as well.  
From \eqref{nexc} we deduce that all non-exceptional entries of 
$\overline{V}_{(n)}$ are uniquely determined from the lower degree terms. 

\begin{lemma}\label{lemma_A}
If $i\neq j$ then the entries of $F_a^{(n)}$ satisfy the following 
equalities: 
\ben
[\overline{F}_a^{(n)}]_{ij}=
[\overline{V}_{(n)}]_{ij}\,\frac{\gl_{a,i}-\gl_{a,j}}{\gl_i-\gl_j}+\ldots,
\een
where the dots stand for terms depending only on $V_{(n')},$ $n'<n$.
\end{lemma}
\proof
We use Taylor's theorem for matrices:
\ben
f(X+Y)=f(X)+d_Xf(Y)+\ldots,
\een
where the dots stand for at least quadratic terms in $Y$. 
Assume now that $X={\rm Diag}(x_1,\ldots,x_N)$ is a diagonal matrix and 
that $f(x)$ is an arbitrary polynomial. Let us compute $d_Xf(Y).$ First, if $f(x)=x^n$ then $d_Xf(Y)=X^{n-1}Y+X^{n-2}YX + \ldots + YX^{n-1}$, i.e.
\beq\label{lemma_A:eq2}
[d_Xf(Y)]_{ij}=
\begin{cases}
Y_{ij}\, \frac{f(x_i)-f(x_j)}{x_i-x_j}, & \mbox{ if } i\neq j \\
Y_{ii}f'(x_i), & \mbox{ if } i=j.
\end{cases}
\eeq
By linearity we get that the above formula holds for all polynomials $f$.

Note that if $f(x)$ is an arbitrary polynomial then 
\beq\label{lemma_A:eq1}
f(V)=f(V_{(0)}) + d_{V_{(0)}}f(V_{(n)})+\ldots,
\eeq
where the dots stand for terms of degree either greater than $n$ or 
terms of degree not exceeding $n$ but depending only on $V_{(n')}$, 
$n'<n.$ On the other hand we have: 
\ben
F_a\d_b = F_a f_b(V_{(0)})e = f_b(V_{(0)})F_a e + [F_a,f_b(V_{(0)})]e,
\een  
 where $e$ is the unity. Note that $F_ae = F_a^{(0)}e.$ Therefore, by
comparing the degree $n$ terms in the above equality and by using 
\eqref{lemma_A:eq1} together with the fact that $f_b(V)$ and $F_a$ 
commute we get: 
\ben
F_a^{(n)}\d_b = [d_{V_{(0)}}f_b(V_{(n)}),F_a^{(0)}]e + \ldots
\een 
where the dots stand for terms depending only on $V_{(n')}$, $n'<n.$ We
multiply both sides of the above equality by $P$ from the left:
\ben
\overline{F}_a^{(n)}\,P\d_b = 
[d_{\overline{V}_{(0)}}f_b(\overline{V}_{(n)}), \overline{F}_a^{(0)}]\,Pe 
+ \ldots.
\een
Since both $\overline{V}_{(0)}$ and $\overline{F}_a^{(0)}$ are diagonal
matrices we can easily get (see \eqref{lemma_A:eq2}) 
\ben
\sum_{s=1}^N[\overline{F}_a^{(n)}]_{is}P_{sb}=
\sum_{\substack{s=1\\s\neq i}}^N
[V_{(n)}]_{is}\, 
\frac{f_b(\gl_i)-f_b(\gl_s)}{\gl_i-\gl_s}\, (\gl_{a,s}-\gl_{a,i})
P_{s1} +\ldots . 
\een
On the other hand we know that $P_{sb}=\gl_{b,s}=f_b(\gl_s)$ and 
$P_{s1}=\gl_{1,s}=1$. Multiply the above equality by $[P^{-1}]_{bj}$ and 
sum over all $b=1,2,\ldots, N:$ 
\ben
[\overline{F}_a^{(n)}]_{ij}= \sum_{\substack{s=1\\s\neq i}}^N
[V_{(n)}]_{is}\, 
\frac{\delta_{i,j}-\delta_{s,j} }{\gl_i-\gl_s}\, (\gl_{a,s}-\gl_{a,i}) 
+\ldots . 
\een
We are given that $i\neq j$ thus $\delta_{i,j}=0$. The only non-zero term
in the above sum is the one corresponding to $s=j$. The lemma follows.
\qed

\begin{lemma}\label{lemma_B}
a) The diagonal entries of $\overline{V}_{(1)}$ are given by 
\ben
[\overline{V}_{(1)}]_{ii}=\sum_{a=1}^N \,
\frac{\gl_{a,i}}{k\gl_i^{k-1}}\,\tau_a.
\een
b) Assume that $(i,j)$ is an exceptional pair of indices. Then  
\ben
[\overline{U}^{(n+1)}]_{ij} = 
\Big(\sum_{a=1}^N \frac{\gl_{a,i}-\gl_{a,j}}{\gl_i-\gl_j}\, \tau_a\Big)
[\overline{V}_{(n)}]_{ij} + \ldots,
\een
where the dots stand for terms depending on $V_{(n')}$, $n'<n$, and the 
non-exceptional entries of $V_{(n)}.$ 
\end{lemma}
\proof
a) We have the following equations: 
$$\overline{U}^{(1)} = 
\overline{V}_{(0)}^{k-1}\overline{V}_{(1)}+ 
\overline{V}_{(0)}^{k-2}\overline{V}_{(1)}\overline{V}_{(0)}+\ldots
+\overline{V}_{(1)}\overline{V}_{(0)}^{k-1}$$ 
and $\d_a\overline{U}^{(1)} = \overline{F}_a^{(0)} + 
[\overline{\mu},\overline{F}_a^{(0)}].$ On the other hand 
$\overline{V}_{(0)}$ and $\overline{F}_a^{(0)}$ are diagonal matrices
and the $i$-th diagonal entries are respectively $\gl_i$ and $\gl_{a,i}$.
Part a) follows. 

b) Assume that $n>1$. Let us compare the $(i,j)$-th, degree $n+1$ entries in the equality $\overline{U}= \overline{V}^k:$
\ben
[\overline{U}^{(n+1)}]_{ij} = \sum_{\substack{0,\leq a_1,a_2,a_3\leq k-2\\
a_1+a_2+a_3 = k-2}} \Big[
\overline{V}_{(0)}^{a_1}\overline{V}_{(1)}\overline{V}_{(0)}^{a_2}
\overline{V}_{(n)}\overline{V}_{(0)}^{a_3} + 
\overline{V}_{(0)}^{a_1}\overline{V}_{(n)}\overline{V}_{(0)}^{a_2}
\overline{V}_{(1)}\overline{V}_{(0)}^{a_3}\Big]_{ij} +\ldots,
\een  
where the dots stand for terms which depend only on $V_{(n')}$, $n'<n.$ Note that if $n=1$ then the second summand should be removed, or equivalently we have to divide the sum by 2. 

The sum is easy to simplify because $\overline{V}_{(0)}$ is a diagonal matrix. We get 
\beq\label{lemma_B:eq1}
\sum_{s=1}^N\, 
\Big(\sum_{\substack{0\leq a_1,a_2,a_3\leq k-2\\
a_1+a_2+a_3 = k-2}} \gl_i^{a_1}\gl_s^{a_2}\gl_j^{a_3} \Big) 
\([\overline{V}_{(1)}]_{is}[\overline{V}_{(n)}]_{sj}+
[\overline{V}_{(n)}]_{is}[\overline{V}_{(1)}]_{sj}\).
\eeq
We want to compute this sum  up to terms independent 
of the exceptional entries of $\overline{V}_{(n)}$. There are three cases.

{\em Case 1:} If $s\neq i$ and $s\neq j$. Then pairs $(s,i)$ and 
$(s,j)$ are either both exceptional or both non-exceptional, because $(i,j)$ is an exceptional pair, (i.e., 
$\gl_i\neq\gl_j$ but $\gl_i^k=\gl_j^k$). We can assume that 
$(s,i)$ and $(s,j)$ are exceptional pairs. In particular,
$\gl_i^k=\gl_j^k=\gl_s^k.$ 
On the other hand, the sum of the $\gl$'s in \eqref{lemma_B:eq1} is 
\ben
\frac{1}{\gl_i-\gl_s}\Big(\frac{\gl_i^k-\gl_j^k}{\gl_i-\gl_j}-
\frac{\gl_s^k-\gl_j^k}{\gl_s-\gl_j}\Big).
\een
The above sum is 0. Hence in case $s\neq i$ and $s\neq j$ there is no
contributions.

{\em Case 2:} $s=j$. Then the sum of the $\gl$'s in \eqref{lemma_B:eq1}
is 
$
-{k\gl_j^{k-1}}/({\gl_i-\gl_j}).
$
On the other hand, since the entry $[V_{(n)}]_{sj}$ is not exceptional and $[V_{(1)}]_{is}$ is already determined we get that the contribution we are interested in is 
\ben
-\frac{k\gl_j^{k-1}}{\gl_i-\gl_j}[V_{(n)}]_{ij}[V_{(1)}]_{jj} = 
-\Big(\sum_{a=1}^N \frac{\gl_{a,j}}{\gl_i-\gl_j}\,\tau_a\Big)[V_{(n)}]_{ij},
\een 
where in the first equality we used part a).

{\em Case 3:} $s=i.$ Just like in the second case we get that the contribution
is 
\ben
-\Big(\sum_{a=1}^N \frac{\gl_{a,i}}{\gl_j-\gl_i}\,\tau_a\Big)[V_{(n)}]_{ij}.
\een
We sum up the contributions from the three cases and then part b) follows. 

If $n=1$, then the contributions from cases 2 and 3  should be doubled, because we have $$[\overline{V}_{(1)}]_{is}[\overline{V}_{(n)}]_{sj}+[\overline{V}_{(n)}]_{is}[\overline{V}_{(1)}]_{sj} = 2[\overline{V}_{(1)}]_{is}[\overline{V}_{(1)}]_{sj}.$$ However this additional factor of 2 is compensated by an earlier division by 2 as it was already explained in the beginning of our proof of part b). 
\qed

Now we are ready to finish the proof of the reconstruction theorem. 
We need to prove that the exceptional entries $[\overline{V}_{(n)}]_{ij}$ are 
uniquely determined in terms of $V_{(n')}$, $n'<n$ and the non-exceptional
entries of $\overline{V}_{(n)}.$ This follows from the equation
\beq\label{exc}
\d_a\overline{U}^{(n+1)} = \overline{F}_a^{(n)} + 
[\overline{\mu},\overline{F}_a^{(n)}]. 
\eeq
Let $(i,j)$ be an exceptional pair and let us compare the $(i,j)$-th
entries in \eqref{exc}. We claim that the $(i,j)$-th entry of
$[\overline{\mu},\overline{F}_a^{(n)}]$ is independent of the 
exceptional entries of $\overline{V}_{(n)}$. Indeed, the entry is given
by
\ben
\sum_{s=1}^n [\overline{\mu}]_{is}[\overline{F}_a^{(n)}]_{sj} - 
[\overline{F}_a^{(n)}]_{is}[\overline{\mu}]_{sj}.
\een
If $s$ equals $i$ or $j$ then $\gl_i^k=\gl_j^k=\gl_s^k$, where the 
first equality holds because $(i,j)$ is an exceptional pair. On the other hand, since $U^{(0)}=V_{(0)}^k$, we have $u^i_{(0)}=\gl_i^k$ for all $i=1,\ldots,N.$ Recalling the second condition of our theorem we get:
$[\overline{\mu}]_{is}=[\overline{\mu}]_{sj}=0.$
If $s$ is different from both $i$ and $j$ then we can assume also that
$\gl_s^k\neq \gl_i^k$ and $\gl_s^k\neq \gl_j^k$, otherwise 
respectively $[\overline{\mu}]_{is}=0$ and $[\overline{\mu}]_{sj}=0.$ In other words $(s,i)$ and $(s,j)$ are not exceptional pairs. According to \leref{lemma_A}, $[\overline{F}_a^{(n)}]_{sj}$ and $[\overline{F}_a^{(n)}]_{is}$ depend only on $V_{(n')}$, $n'<n$ and the non-exceptional entries of $\overline{V}_{(n)}$. 

To finish the proof it remains only to recall \leref{lemma_A}, 
part b) of \leref{lemma_B}, and \eqref{exc}. We get 
\ben
\Big(\sum_{a=1}^N \frac{\gl_{a,i}-\gl_{a,j}}{\gl_i-\gl_j}\,\tau_a
\Big)\,
\d_{\tau_a} [V_{(n)}]_{ij} = \mbox{ known terms }.
\een    
Notice that the above sum is non-zero because the coefficient in front of $\tau_2$ is 1. 
\qed

\subsection{Twisted curves} 
In this Section every scheme is over $\com$, and the terms ``orbifold'' and ``smooth Deligne-Mumford stack'' are used interchangeably. Let $C$ be a smooth curve, $p_1,...,p_n\in C$ distinct points, and $k_1, ...., k_n$ positive integers. Given these data, we consider the stack $C[(p_1,k_1),...,(p_n,k_n)]$, which is constructed as the stack of roots of line bundles on $C$. More precisely, the stack $C[(p_1,k_1),...,(p_n,k_n)]$ is the fiber product $\radice{k_1}{(\sO_{p_1}, \sigma_1)/C}\times_C...\times_C \radice{k_n}{(\sO_{p_n},\sigma_n)/C}$. The stack $\radice{k_i}{(\sO_{p_i},\sigma_i)/C}$ is the stack of $k_i$-th root of the line bundle $\sO(p_i)$ with the canonical section $\sigma_i:\sO\to \sO(p_i)$ . An object of $\radice{k_i}{(\sO_{p_i}, \sigma_i)/C}$ over a $C$-scheme $T$ is 
\begin{itemize}
\renewcommand{\labelitemi}{--}
\item
a line bundle $M$;
\item
an isomorphism $\phi$ of $M^{\otimes k_i}$ with the pullback of $\sO(p_i)$ via $T\to C$;
\item
a section $\tau$ of $M$ such that $\phi(\tau^{k_i})=\sigma_i$.
\end{itemize}

More details of this construction can be found in \cite{agv2} and \cite{ca}. An alternative description of $C[(p_1,k_1),...,(p_n,k_n)]$ using log geometry may be found in \cite{ol}. \'Etale locally near a point $p\in C\setminus \{p_1,...,p_n\}$, the stack $C[(p_1,k_1),...,(p_n,k_n)]$ is isomorphic to the curve $C$. \'Etale locally near the point $p_i$, the stack $C[(p_1,k_1),...,(p_n,k_n)]$ is isomorphic the the stack quotient $[\text{Spec}\,\mathbb{C}[x]/\mu_{k_i}]$ where the group $\mu_{k_i}$ acts via $x\mapsto \zeta x$ for $\zeta\in \mu_{k_i}$. The natural projection $$C[(p_1,k_1),...,(p_n,k_n)]\to C$$ exhibits $C$ as its coarse moduli space.

\subsection{Orbifold quantum cohomology}
Our focus is a simple case of this construction, namely $$\mC_{k,m}:=\proj^1[(0,k),(\infty,m)],$$ for integers $k, m\geq 1$. We call $\mC_{k,m}$ a $2$-pointed $\proj^1$-orbifold. Roughly speaking, this is the curve $\proj^1$ with orbifold points $B\mathbb{Z}_k$ and $B\mathbb{Z}_m$ at $0$ and $\infty$ respectively\footnote{Of course placing the two orbifold points elsewhere on $\proj^1$ results isomorphic orbifolds.}. Note that for $k, m$ coprime, the orbifold $\proj^1[(0,k),(\infty,m)]$ is isomorphic to the weighted projective line $\proj^1(k,m)$. If $k$ and  $m$ are not coprime, then the weighted projective line $\proj^1(k,m)$ has nontrivial generic stabilizer. But  $\proj^1[(0,k),(\infty,m)]$ always has trivial generic stabilizers, thus it is not isomorphic to $\proj^1(k,m)$. Also, it is obvious that $\proj^1[(0,k),(\infty,m)]\simeq \proj^1[(0,m),(\infty,k)]$.

To each orbifold one can associate another orbifold, the so-called {\em inertia orbifold}, which plays a key role in orbifold Gromov-Witten theory. By definition, the inertia orbifold of a given orbifold $\X$ is defined to be the fiber product (in the $2$-category of stacks) $I\X:=\X\times_{\Delta,\X\times \X,\Delta} \X$, where $\Delta: \X\to \X\times \X$ is the diagonal morphism. In categorical terms, the objects of $I\X$ are:
$$Ob(I\X):=\{(x,g)|x\in Ob(\X), g\in Aut(x)\}$$
$$=\{(x,H,g)|x\in Ob(\X), H\subset Aut(x), g \mbox{ a generator of } H\}.$$
There are two natural maps: $q: I\X\to \X$ given by forgetting the choice of $g\in Aut(x)$, and $I: I\X\to I\X$ given by $g\mapsto g^{-1}$.
The inertia orbifold $I\X$ is disconnected (unless $\X$ is a connected manifold). Write $I\X=\coprod_{i\in \sI} \X_i$ for the decomposition into connected components. On each connected component $\X_i$ there is a trivial action of the cyclic group $\mathbb{Z}_{r_i}$. Thus $\mathbb{Z}_{r_i}$ acts on the vector bundle $q^*T_\X$. This yields a decomposition $q^*T_\X=\oplus_j E_{ij}$ into eigen-bundles, where $\mathbb{Z}_{r_i}$ acts on $E_{ij}$ via multiplication by $\exp(2\pi\sqrt{-1}\frac{k_j}{r_i})$ with $0\leq k_j<r_i$. The age associated to the component $\X_i$ is defined to be $age(\X_i):=\sum_j k_j/r_i$.

In our case, it is easy to see that $$I\mC_{k,m}\simeq\mC_{k,m}\cup\bigcup_{1\leq i\leq k-1} B\mu_k(i)\cup\bigcup_{1\leq j\leq m-1}B\mu_m(j).$$
Here for each $i, j$ we have $B\mu_k(i)\simeq B\mu_k$ and $B\mu_m(j)\simeq B\mu_m$. The age associated to the component $\mC_{k,m}$ is $0$, the age associated to $B\mu_k(i)$ is $i/k$, the age associated to $B\mu_m(j)$ is $j/m$.

Next we turn to orbifold cohomology. As a graded vector space, the orbifold cohomology of $\X$ (with complex coefficients) is defined to be $H_{orb}^*(\X):=H^*(I\X)$ with the grading defined as follows: a class $a\in H^p(\X_i)$ is assigned the degree $p+2age(\X_i)$. The orbifold cohomology of $\mC_{k,m}$, as a vector space, is given by $$H_{orb}^*(\mC_{k,m})=H^0(\mC_{k,m})\oplus H^2(\mC_{k,m})\oplus \oplus_{1\leq i\leq k-1}H^0(B\mu_k(i))\oplus \oplus_{1\leq j\leq m-1}H^0(B\mu_m(j)).$$
An element in $H^0(B\mu_k(i))$ is assigned degree $2i/k$, and an element in $H^0(B\mu_m(j))$ is assigned degree $2j/m$. This gives $H_{orb}^*(\mC_{k,m})$ the structure of a graded vector space.

We fix some notations.  Let $1\in H^0(\mC_{k,m})$ be the Poincar\'e dual of the fundamental class, $p\in H^2(\mC_{k,m})$ the Poincar\'e dual of a point, $x_i\in H^0(B\mu_k(i))$ the Poincar\'e dual of the fundamental class for each $i$, and $y_i\in H^0(B\mu_m(j))$ the Poincar\'e dual of the fundamental class for each $j$. 

In general the orbifold cohomology space $H_{orb}^*(\X)$ carries a non-degenerate pairing $(\,\,,\,\,)_{orb}$ called orbifold Poincar\'e pairing. It is defined as follows: for $a,b\in H^*(I\X)$, define $(a,b)_{orb}:=\int_{I\X}a\wedge I^*b$. In our case this pairing is given as follows: 
$$(x_{i},x_{k-i})_{orb}=1/k, (y_{j},y_{m-j})_{orb}=1/m, (1, p)_{orb}=1=(p,1)_{orb}; \text{ and }0\text{ otherwise}.$$

A recent advance in the study of orbifolds is that orbifold cohomology $H_{orb}^*(\X)$ carries a nontrivial ring structure called orbifold cup product, see \cite{cr1} and \cite{agv1}. We briefly recall its definition. The geometric object central to the construction of this ring structure, as well as orbifold Gromov-Witten theory, is the notion of orbifold stable maps. An orbifold stable map $f: \mathfrak{C}\to \X$ is a representable map from a nodal curve $\mathfrak{C}$, possibly having orbifold structures at marked points and nodes, to the orbifold $\X$. We may fix discrete invariants and consider moduli spaces\footnote{A technical point: we consider here orbifold stable maps with sections to all gerbes. Our notation here agrees with that in \cite{agv1} and \cite{ts0}.} $\Mbar_{g,n}(\X,d)$ parametrizing $n$-pointed orbifold stable maps of genus $g$ and degree $d\in H_2(X,\mathbb{Z})$. These moduli spaces come with two kinds of maps: the evaluation map at the $i$-th marked point $ev_i: \Mbar_{g,n}(\X,d)\to I\X$; the map $\pi:\Mbar_{g,n}(\X,d)\to \Mbar_{g,n}(X,d)$ given by passing to coarse moduli spaces.

Deformation theory of orbifold stable maps yields a perfect obstruction theory on $\Mbar_{g,n}(\X,d)$, from which one can construct a virtual fundamental class $[\Mbar_{g,n}(\X,d)]^{vir}\in H_*(\Mbar_{g,n}(\X,d),\mathbb{Q})$, see \cite{agv2}.

Now we can define orbifold cup products: for $a,b\in H_{orb}^*(\X)$, define $$a\cdot b:=(I\circ ev_3)_*(ev_1^*a\cup ev_2^*b\cap [\Mbar_{0,3}(\X,0)]^{vir}).$$ This product respects gradings, making $(H_{orb}^*(\X), \cdot)$ a graded commutative associative $\com$-algebra.

We now describe the orbifold cohomology ring structure of $H_{orb}^*(\mC_{k,m})$ as explained above. By definition, classes in $H^*(\mC_{k,m})=H^*(\proj^1)$ multiply as usual. Furthermore, using only the definition, we have $$x_i\cdot y_j=0 \text{ for every }i, j;$$
$$x_{i_1}\cdot x_{i_2}=x_{i_1+i_2} \text{ if }i_1+i_2\leq k-1;$$
$$y_{j_1}\cdot y_{j_2}=y_{j_1+j_2} \text{ if }j_1+j_2\leq m-1;$$
$$k{x_1}^k=m{y_1}^m=p.$$
It follows that, as rings, $$H_{orb}^*(\mC_{k,m})\simeq \mathbb{C}[x,y]/(kx^k-my^m, xy),$$ where we identify $x_1=x$ and $y_1=y$.

\begin{remark}{\em
The calculation of orbifold Poincar\'e pairing and orbifold cup product for $\mC_{k,m}$ can be easily generalized to the more general twisted curve $C[(p_1,k_1),...,(p_n,k_n)]$. We won't need this here.}
\end{remark}

The definition of orbifold cup product involves only degree $0$ orbifold stable maps. Intersection numbers on moduli spaces of orbifold stable maps of nonzero degrees can be packaged to give a deformation of the orbifold cup product, which we now describe. 

For classes $a_1,...,a_n\in H_{orb}^*(\X)$, define the {\em genus zero primary orbifold Gromov-Witten invariant} $\<a_1,..,a_n\>_{0,n,d}$ to be the integral $$\int_{[\Mbar_{0,n}(\X,d)]^{vir}}ev_1^*a_1\wedge...\wedge ev_n^*a_n.$$ Fix an additive basis $\{\phi_\alpha\}$ of $H_{orb}^*(\X)$ and write $\{\phi^\alpha\}$ for its dual basis. For classes $a,b\in H_{orb}^*(I\X)$, the formula 
$$a\ast_t b:=\sum_{n,d}\frac{Q^d}{n!}\<a,b,\phi_\alpha,t,...,t\>_{0,n+3,d}\phi^\alpha,$$
defines a ring structure on $H_{orb}^*(\X)$ with coefficient ring enlarged to the Novikov ring $\com[[H_2(\X)]]$. This product $\ast_t$ depends on a parameter $t\in H_{orb}^*(\X)$. Associativity of $\ast_t$ is nontrivial. The ring $BQH_{orb}^*(\X):=(H_{orb}^*(\X),\ast_t)$ is called the big orbifold quantum cohomology ring. The variables $Q$ in the Novikov ring are assigned degrees so that $deg(Q^d)=2\int_d c_1(T_\X)$. The product $\ast_t$ respects degrees.

When restricting to $t\in H^2(\X)$, an easy application of the divisor equation shows that $\ast_t$ may be identified with the following product:
$$a \star_t b:=\sum_d (Qe^t)^d (I\circ ev_3)_*(ev_1^*a\cup ev_2^*b\cap [\Mbar_{0,3}(\X,d)]^{vir}).$$ We call $QH_{orb}^*(\X):=(H_{orb}^*(\X),\star_t)$ the small orbifold quantum cohomology ring.

We now describe the small orbifold quantum cohomology ring $QH_{orb}^*(\mC_{k,m})$. First note that the Picard group $\text{Pic}(\mC_{k,m})$ is generated by two line bundles $L_0, L_\infty$ such that $L_0^{\otimes k}\simeq L_\infty^{\otimes m}$ and both are isomorphic to the pull-back of $\sO_{\mathbb{P}^1}(1)$. We have $deg\, L_0=1/k, deg\, L_\infty=1/m$. The canonical line bundle $K_{\mC_{k,m}}\simeq L_0^\vee\otimes L_\infty^\vee$ has degree $-1/k-1/m$. Additively we have $QH_{orb}^*(\mC_{k,m})=H_{orb}^*(\mC_{k,m})\otimes_\mathbb{C}\mathbb{C}[[q]]$, where the variable $q:=Qe^{t}, t\in H^2(\mC_{k,m})$ is assigned degree $2/k+2/m$. The product structure of $QH_{orb}^*(\mC_{k,m})$ is a deformation of that on $H_{orb}^*(\mC_{k,m})$. So we only have to analyze how to deform the relations in $H_{orb}^*(\mC_{k,m})$. By degree consideration, the relation $xy=0$ can only be deformed to $xy=cq$. Here by definition $c$ is the orbifold Gromov-Witten invariant $\<x,y, p\>_{0, 3, 1}$, which is clearly equal to $1$. 

The relation $kx^k-my^m=0$ remains unchanged. This is easily seen by degree consideration if $k,m$ are co-prime. In general, it follows from the fact that the {\em classes} $kx^k, my^m$ in $H_{orb}^*(\mC_{k,m})$ are both equal to $p$, and the following 

\begin{lemma}
\hfill
\begin{enumerate}
\item
The product $x^{\star_t a}$ of $a$ copies of $x$ is equal to $x^a$ if $1\leq a\leq k$.
\item
The product $y^{\star_t b}$ of $b$ copies of $y$ is equal to $y^b$ if $1\leq b\leq m$.
\end{enumerate}
\end{lemma}

\begin{proof}
We only prove the statement about $x^{\star_t a}$, an analogous argument proves the statement about $y^{\star_t b}$. 

We need the following non-vanishing conditions:
\begin{equation}\label{non-vanishing_1}
\text{If } \<x^i, x^j, x^l\>_{0,3, d}\neq 0, \text{ then } i+j+l\equiv d\, (\text{mod } k), \text{ and } d\equiv 0\, (\text{mod }m).
\end{equation}

\begin{equation}\label{non-vanishing_2}
\text{If } \<x^i, x^j, y^l\>_{0,3,d}\neq 0, \text{ then } i+j\equiv d\, (\text{mod } k), \text { and } d\equiv l\, (\text{mod }m).
\end{equation}

We first prove the statement about $x^{\star_t a}$ assuming (\ref{non-vanishing_1}), (\ref{non-vanishing_2}). We proceed by induction on $a$. Clearly $x^{\star_t 1}=x$. Suppose that $x^{\star_t a}=x^a$ for some $a\leq k-1$. We may write 
$$x^{\star_t a+1}=x^a\star_t x=x^{a+1}+\sum_{d > 0}q^d \left(a_{0, d}1+b_{0,d}p+\sum_{i=1}^{k-1}a_{i,d}x^i+\sum_{j=1}^{m-1}b_{j,d}y^j \right).$$
Note that the right side of the equation above should be homogenous of degree $2(a+1)/k$.

Suppose that $a_{i,d}\neq 0$ for some $1\leq i\leq k-1$. Then by definition we have $\<x^a,x,x^{k-i}\>_{0,3,d}\neq 0$.  By (\ref{non-vanishing_1}), we have $a+1\equiv i+d\, (\text{mod }k)$ and $d\equiv 0\, (\text{mod }m)$. In particular $d/m\geq 1$. Now by comparing degrees, we find 
$$1\geq \frac{a+1}{k}=\frac{i}{k}+d \left(\frac{1}{k}+\frac{1}{m}\right)=\frac{i+d}{k}+\frac{d}{m} >1,$$ which is a contradiction. The same argument proves that $a_{0,d}=0$.

Suppose that $b_{j,d}\neq 0$ for some $1\leq j\leq m-1$. Then by definition we have $\<x^a,x,y^{m-j}\>_{0,3,d}\neq 0$. By (\ref{non-vanishing_2}), we have $a+1\equiv d\,(\text{mod } k)$ and $m-j\equiv d\, (\text{mod }m)$. In particular $(d+j)/m\geq 1$. Again by comparing degrees, we find 
$$1\geq \frac{a+1}{k}=\frac{j}{m}+d\left(\frac{1}{k}+\frac{1}{m}\right)=\frac{d}{k}+\frac{d+j}{m}>1,$$ which is a contradiction.

Finally, $b_{0,d}=0$ because the degree of $p\cdot q^d$ is $2(1+d(1/k+1/m))>2$, while the degree of $x^{a}\star_t x$ is at most $2$.

Now we prove (\ref{non-vanishing_1}). If $ \<x^i, x^j, x^l\>_{0,3, d}\neq 0$, then the relevant moduli space must be non-empty. So there exists a three-pointed, degree $d$ orbifold stable map $f: \mathfrak{C}\to \mC_{k,m}$ with stack structures on $\mathfrak{C}$ prescribed by the insertions. The holomorphic Euler characteristics $\chi (\mathfrak{C}, f^*L_0)$ and $\chi(\mathfrak{C}, f^*L_\infty)$ are integers. By Riemann-Roch, we find
\begin{equation*}
\begin{split}
&\chi(\mathfrak{C}, f^*L_0)=1+\frac{d}{k}-\frac{i}{k}-\frac{j}{k}-\frac{l}{k},\\
&\chi(\mathfrak{C}, f^*L_\infty)=1+\frac{d}{m}.
\end{split}
\end{equation*}
The result follows.

The proof of (\ref{non-vanishing_2}) is similar: if $ \<x^i, x^j, y^l\>_{0,3,d}\neq 0$, then the relevant moduli space is not empty. So there exists a three-pointed, degree $d$ orbifold stable map $f: \mathfrak{C}\to \mC_{k,m}$ with stack structures on $\mathfrak{C}$ prescribed by the insertions. One calculates by Riemann-Roch that, in this case, 
\begin{equation*}
\begin{split}
&\chi(\mathfrak{C}, f^*L_0)=1+\frac{d}{k}-\frac{i}{k}-\frac{j}{k},\\
&\chi(\mathfrak{C}, f^*L_\infty)=1+\frac{d}{m}-\frac{l}{m}.
\end{split}
\end{equation*}
The result follows by integrality of $\chi(\mathfrak{C}, f^*L_0)$ and $\chi(\mathfrak{C}, f^*L_\infty)$.

\end{proof}

Hence we obtain the following presentation of the small orbifold quantum cohomology ring:
\begin{equation}\label{smallqh}
QH_{orb}^*(\mC_{k,m})\simeq \mathbb{C}[[q]][x,y]/(kx^k-my^m, xy-q).
\end{equation}

This presentation allows us to set $Q$ to any nonzero complex number. We do so from now on.

Note that in case of $k,m$ coprime, (\ref{smallqh}) coincides with the calculations in \cite{agv2} for weighted projective lines $\proj^1(k,m)$.

\subsection{Frobenius structure}
It is known that genus zero Gromov-Witten theory provides a natural Frobenius structure on the cohomology of the target space. The same is true for orbifolds: orbifold cohomology $H_{orb}^*(\X)$ of an orbifold $\X$ carries a natural Frobenius structure arising from genus zero orbifold Gromov-Witten invariants. The ingredients of this Frobenius structure are summarized as follows.
\begin{itemize}
\item
The space on which the Frobenius structure is based: the orbifold cohomology $H_{orb}^*(\X)\otimes_\mathbb{C}\com[[H_2(\X)]]$;
\item
the flat metric is given by the orbifold Poincar\'e pairing $(\,\,,\,\,)_{orb}$; 
\item
the product structure is given by the orbifold big quantum product $\star_t$. 
\end{itemize}

We now turn to the special case $BQH_{orb}^*(\mC_{k,m})$. Consider the following 
homogeneous additive basis of $H_{orb}^*(\mC_{k,m})$,
\begin{equation}\label{flatbasis}
 x_{k-1},...,x_{1},1 , y_{1},...,y_{m-1}, p.
\end{equation}
Note the ordering of these classes. In this basis, we may write a class in 
$H_{orb}^*(\mC_{k,m})$ as $$\sum_{i=1}^{k-1}s_ix_{k-i}+s_k1+
\sum_{j=1}^{m-1}s_{k+j}y_{j}+s_Np,$$ 
where $N=k+m$. 

By expressing multiplications by $x$ and $y$ in the presentation (\ref{smallqh}) as matrices using the (ordered) basis (\ref{flatbasis}), it is easy to show that the Frobenius manifold $BQH_{orb}^*(\mC_{k,m})$ is semi-simple along $H^2(\mC_{k,m})$.
 
The coordinates $\{s_1,..., s_N\}$ are flat coordinates of $BQH_{orb}^*(\mC_{k,m})$. 
In these coordinates, the Euler vector field reads
$$s_k\partial_{s_k}+\sum_{i=1}^{k-1}(i/k)s_i\partial_{s_i}+
\sum_{j=1}^{m-1}(1-j/m)s_{k+j}\partial_{s_{k+j}}+(1/k+1/m)\partial_{s_N}.$$

{\em Proof of Theorem \ref{t3}.} 
Using Theorem \ref{reconstruction}, we prove that the following map
\begin{equation}\label{isomap}
s_i\mapsto \tau_i, \text{ for } i=1, ...., N-1, s_N\mapsto \tau_N/m,
\end{equation} 
is an isomorphism between the Frobenius structures respectively on the big quantum cohomology $BQH_{orb}^*(\mC_{k,m})$ and on $M_{k,m}$.

The map (\ref{isomap}) identifies $M_{k,m}$ and $BQH^*_{\rm orb}(\CC_{k,m})$
as complex manifolds. It also identifies the corresponding flat metrics,
unity vector fields and Euler vector fields. It remains only to verify that both Frobenius structures satisfy the conditions of \thref{reconstruction}. 

Note that at the point $\tau_1=\ldots =\tau_N = 0$ we have isomorphisms of Frobenius algebras $T_0 M_{k,m}\iso T_0BQH^*_{orb}(\CC_{k,m})\iso \C[x,y]/\< kx^k - my^m, xy-Q\>,$ by (\ref{smallqh}). Up to a scalar, the $k$-th root of the restriction of the Euler vector field to $T_0M_{k,m}$ is given by $x$. From the above presentation of $T_0M_{k,m}$, it follows that $x$ is an invertible generator of the Frobenius algebra and that the point $\tau=0$ is semisimple, i.e., the first condition in \thref{reconstruction} is satisfied. It remains only to verify the second one. 

Pick the following basis of $T_0M_{k,m}$:
\ben
\phi_1=x^{k-1},\ldots,\phi_k=1,\phi_{k+1}=\frac{Q}{x},\ldots,\phi_{k+m}=\Big(\frac{Q}{x}\Big)^m.
\een 
It is easy to see that the eigenvalues of $x\bullet_{0}$ are given by:
\ben
\gl_i = \Big( mQ^m/k\Big)^{1/k} \exp\(2\pi\sqrt{-1}\,i/N\),\quad 1\leq i\leq N.
\een
On the other hand, the quantum product $\bullet_0$ is diagonalized by the matrix $P$ whose $i$-th row is given by 
$\Big(\gl_i^{k-1},\ldots,1,Q/\gl_i,\ldots,(Q/\gl_i)^m\Big).$ 
The $i$-th column of the matrix inverse to $P$ is given by:
$\frac{1}{N}\Big(\gl_i^{-(k-1)},\ldots,1,\gl_i/Q,\ldots,(\gl_i/Q)^m\Big)$ (recall that by definition $N=k+m$). We need to prove that if $i$ and $j$ are such that $\gl_i^k=\gl_j^k$ then the $(i,j)$-th entry of $\overline{\mu}:=P\mu P^{-1}$ is 0. On the other hand in the above basis $\mu$ is represented by the diagonal matrix 
\ben
{\rm Diag}\Big(i/k-1/2,1\leq i\leq k-1,\ 1/2-j/m, 0\leq j\leq m\Big).
\een
Therefore we have to verify that 
\beq\label{spectral_identity}
\sum_{s=1}^{k-1} P_{is}\(\frac{s}{k}-\frac{1}{2}\)[P^{-1}]_{sj} + 
\sum_{s=0}^{m} P_{i,k+s}\(\frac{1}{2}-\frac{s}{m} \)[P^{-1}]_{k+s,j} =0.
\eeq
Note that $P_{is}[P^{-1}]_{sj}=(\gl_i/\gl_j)^{k-s}/N = (\gl_j/\gl_i)^s/N,$ where for the second equality we used that $\gl_i^k=\gl_j^k.$ Put $x=\gl_j/\gl_i$. Using that $x^k=1$ and the identity:
\ben
x+2x+\ldots (k-1)x^{k-1} = \frac{(k-1)x^{k+1}-kx^k + x}{(x-1)^2},
\een 
we get that the first sum in \eqref{spectral_identity} equals $\frac{1}{N}\(\frac{1}{2} + \frac{1}{x-1}\).$ For the second sum we have 
$P_{i,k+s}[P^{-1}]_{k+s,j}=(\gl_j/\gl_i)^s/N.$ Note that the summands corresponding to $s=0$ and $s=m$ cancel each other, so we may assume that the summation range is from 1 to $m-1$. Also, we have that $x^m=1,$ because $N=k+m,$ $x^k=1$ and $x^N=1.$ So the second sum simplifies to $-\frac{1}{N}\(\frac{1}{2} + \frac{1}{x-1}\).$

\qed

\subsection{Descendent potential}
We recall the definition of the {\em descendent orbifold Gromov-Witten invariants}, 
which plays an important role in orbifold Gromov-Witten theory.

Recall that on the moduli space $\Mbar_{g,n}(X,d)$ of stable maps to the {\em coarse moduli space} $X$ there are $n$ line bundles $L_1,...,L_n$ associated to the marked points. The fiber of $L_i$ at a moduli point $(f: (C,p_1,...,p_n)\to X)$ is the cotangent space $T_{p_i}^*C$. Consider the pullback line bundles $\pi^*L_i$. The descendent classes in orbifold Gromov-Witten theory of $\X$ are defined to be $\bpsi_i:=c_1(\pi^*L_i)$.

The totality of descendent orbifold Gromov-Witten invariants can be packaged in a generating function, called the total descendent potential of $\X$, which is defined as follows:
$$\D_\X:=\exp\left(\sum_{g\geq 0} \epsilon^{2g-2}\sum_{n,d}\frac{Q^d}{n!}\int_{[\Mbar_{g,n}(\X,d)]^{vir}}\bigwedge_{i=1}^n \sum_{k=0}^\infty ev_i^*t_k\bpsi_i^k\right).$$
The total descendent potential $\D_\X$ is viewed as a function of $\bt:=\sum_{k\geq 0}t_kz^k$ and, via the dilaton shift $\bq=\bt-1z$, as an element in the Fock space--the space of functions on $H_{orb}^*(\X)[z]$. 
Assume that $\X$ has semi-simple orbifold quantum cohomology and denote by 
$\D^{{\rm BQH}(\X)}$ the descendent potential corresponding to the semi-simple
Frobenius structure (\cite{G1}). Then we have the orbifold version of Givental's conjectural formula: $\D_\X=\D^{{\rm BQH}(\X)}.$ Recent work of C. Teleman \cite{te} is very close to providing a proof of this. In our case  (i.e., $\X=\mC_{k,m}$) this formula can be proven by virtual localization. Details will be given in \cite{ts}.

In conclusion, we formulate the following conjecture:
\begin{conjecture}\label{star}
{\em  The total descendent potential of $\mC_{k,m}$ is a tau-function of the extended 
bi-graded Toda hierarchy corresponding to $M_{k,m}$ (\cite{car}).}
\end{conjecture}

Once the HQE \eqref{HQE:descendents} is shown to describe the extended bigraded Toda hierarchy, Conjecture \ref{star} will follow from results in this paper and Givental's formula.

We remark that the extended bigraded Toda hierarchy with $k=m=1$ coincides with the extended Toda hierarchy. However, the HQE \eqref{HQE:descendents} specialized to $k=m=1$ are different from the HQE for extended Toda hierarchy given in \cite{M2}. Therefore, the case $k=m=1$ of Conjecture \ref{star} provides yet another 
formulation of the Toda conjecture about Gromov--Witten invariants of 
$\C \proj^1.$


\section{Period vectors and vertex operators}\label{vertexoper}


\subsection{Period vectors near a critical value}
Let $\tau\in M$ be a semi-simple point, i.e., $f_\tau$ 
is a Morse function. A relative cycle 
$\gb\in H_1(\C^*,V_{\tau,\gl_0};\Z)$ is called a Lefschetz thimble corresponding 
to a path $C_i$ from $\gl_0$ to a critical value $u_i$ of $f_\tau$ if $\gb$ is 
represented by the two components in $f_\tau^{-1}(C_i)$ which meet at the  
critical point above $u_i.$ 

\begin{lemma}\label{laplace_transform}
Let $C_\infty$ be a path from $\gl_0$ to $\gl=-z\cdot(+\infty)$  and denote by 
\ben
\B\in \lim_{M\rightarrow\infty} H_1(\C^*,\{x\in \C^*:{\rm Re}( f_\tau/z)<-M\};\Z)
\een
the cycle obtained from $\gb$ by a parallel transport along $C_\infty.$ 
Then 
\ben
J_\B(\tau,z) = (-2\pi z)^{-1/2}\int_{u_i}^{-z\cdot(+\infty)} e^{\gl/z}
I_\gb^{(0)}(\tau,\gl)d\gl.
\een 
\end{lemma}
\proof
The oscillating integral $(-2\pi z)^{1/2}\J_\B$ can be transformed
as follows.
\ben
\int_\B e^{f_\tau/z}\omega & = & \int_{u_i}^{-z\cdot(+\infty)} e^{\gl/z}
\( \int_{\d\gb}\frac{\omega}{df_\tau}\) d\gl = 
\int_{u_i}^{-z\cdot(+\infty)} e^{\gl/z}\d_\gl
\( \int_{\d\gb} d^{-1}\omega\) d\gl \\
& = &
-z^{-1} \int_{u_i}^{-z\cdot(+\infty)} e^{\gl/z}\( \int_{\gb}\omega\) d\gl,
\een
where in the last equality we applied integration by parts and 
the Stokes' formula. The lemma follows because, by definition,
\ben
(J_\B,\d_i)= z\d_i\J_\B\quad \mbox{ and } \quad
(I_\gb^{(0)},\d_i) = -\d_i\int_{\gb}\omega.
\een
\qed

\begin{lemma}\label{I_ui}
Let $\xi$ be sufficiently close to $u_i$ Then
\beq\label{i_0}
I^{(0)}_{\gb}(\tau,\xi) = 
\frac{2}{\sqrt{2(\xi-u_i)}}\({\bf 1}_i + A_{i,1}[2(\xi-u_i)]+ 
A_{i,2}[2(\xi-u_i)]^2+\ldots\ \), 
\eeq
where the path $C_i'$ specifying $\gb(\tau,\xi)$ is the same as $C_i$ except for
the end where the two paths split: $C_i$ leads to $u_i$ and $C_i'$ leads to 
$\xi.$  
\end{lemma} 
\proof
We follow \cite{AGV}, chapter 3, section 12, Lemma 2. 
In a neighborhood of the critical point above $u_i,$ we choose a {\em unimodular} 
coordinate $y$ for the volume form $\omega$ i.e., $\omega = dy.$ The Taylor's expansion of $f_\tau$ is
\ben
f_\tau(y) = u_i + \frac{\Delta_i}{2}(y-y_i)^2 + \ldots ,
\een
where $y_i$ is the $y$-coordinate of the critical point $q_i$ corresponding 
to $u_i.$ From this
expansion we find that the equation $f_\tau(y) = \xi$ has two 
solutions in a neighborhood of $y=q_i$:
\ben
y_\pm = y_i \pm \frac{1}{\sqrt\Delta_i} \sqrt{2(\xi-u_i)} + \mbox{ h.o.t.},
\een
where h.o.t. means higher order terms. Thus the integral in the definition
of $I_{\gb}^{(0)}(\tau,\xi)$ has the following expansion:
\ben
\int_{\gb(\tau,\xi)}\omega = y_+(\tau,\xi) - y_-(\tau,\xi) = 
\frac{2}{\sqrt\Delta_i} \sqrt{2(\xi-u_i)} + \mbox{ h.o.t. },
\een
which yields
\ben
(I_{\gb}^{(0)}(\tau,\xi),\d_j) = 
\frac{2}{\sqrt{2(\xi-u_i)}}\, 
\frac{1}{\sqrt\Delta_i}\d_j u_i + \mbox{ h.o.t. }.
\een
The lemma follows because
\ben
{\bf 1}_i =\sqrt{\Delta_i}\d_{u_i}= \frac{1}{\sqrt\Delta_i}du_i = 
\sum_{j=1}^N \frac{1}{\sqrt\Delta_i}(\d_j u_i)d\tau_j.
\een
\qed

\begin{lemma}\label{vector}
Let $\gl$ be close to $u_i$. Then 
\beq\label{psi_r}
\f_\tau^{\gb/2}(\gl) =\Psi R \ 
\sum_{n\in \Z} (-z\d_\gl)^n \frac{e_i}{\sqrt{2(\gl-u_i)}}, 
\eeq
where $e_i\in \C^N\iso T_\tau M$ corresponds to ${\bf 1}_i$.
\end{lemma} 
\proof
Using \leref{laplace_transform}, we will compute the stationary phase
asymptotic of $J_\B.$ The computation is the same as in the proof of 
Theorem 3 in \cite{G2}. 

Near the critical value the period  $I_{\gb}^{(0)}(\tau,\gl)$ 
has the expansion \eqref{i_0}. Using the change of variables 
$\gl-u_i = -zx^2/2$ we compute
\ben
\frac{2}{\sqrt{-2\pi z}}
\int_{u_i}^{-z\cdot(+\infty)} e^{\gl/z}
[2(\gl-u_i)]^{k-1/2} d\gl & = & (-z)^k\frac{e^{u_i/z}}{\sqrt{2\pi}}
\int_{-\infty}^\infty e^{-x^2/2}x^{2k}dx \\
&=& e^{u_i/z}(-z)^k (2k-1)!! \ . 
\een  
Thus $J_\B$ has the following asymptotic:
\ben
J_\B \sim \(\sum_{k=0}^\infty (2k-1)!!\, A_{i,k} (-z)^k\)e^{u_i/z}. 
\een
Since by definition the asymptotic of $J_\B$ is $\Psi R e^{U/z}e_i$, we get 
$(-1)^k(2k-1)!!A_{i,k}={\Psi R_k}e_i $. Thus
{\allowdisplaybreaks
\ben
\f_\tau^{\gb}(\gl) & = & 
\sum_{n\in\Z}(-z\d_\gl)^n I^{(0)}_{\gb}(\tau,\gl)  \\
& = &
2\sum_{n\in \Z}\sum_{k=0}^\infty 
(-z\d_\gl)^n (-1)^k\Psi R_k{e}_i \frac{[2(\gl-u_i)]}{(2k-1)!!}^{k-1/2} \\
&=&
2\sum_{n\in\Z}\sum_{k=0}^\infty (-z\d_\gl)^n  \Psi R_k
(-\d_\gl)^{-k} \frac{e_i}{\sqrt{2(\gl-u_i)}}  \\
&=& 
2 \(\Psi R \) \sum_{n\in\Z}(-z\d_\gl)^n\frac{e_i}
{\sqrt{2(\gl-u_i)}}.
\een}
\qed


\subsection{Differential equations}

Let $(\tau,\gl)\in (M\times\C)\backslash\Delta.$ We will show that 
$\f_\tau^*(\gl),$ where $*$ is either a relative cycle or a one-point 
cycle, satisfies a certain system of differential equations which, in 
some sense, uniquely determines the corresponding periods.  

Note that 
\ben
\frac{1}{n!} d^{-1}\((\gl-f_\tau)^n \omega \) = \Phi^{(n)}_{\tau,\gl}
+\frac{1}{n!} \( \int_\phi \(\gl-f_\tau\)^n \omega\)\log x, 
\een
where $\Phi^{(n)}_{\tau,\gl}$ is a Laurent polynomial in $x$. Moreover, $\Phi^{(n)}_{\tau,\gl}$ is homogeneous of degree $n$  with respect to our grading conventions.
\begin{lemma}\label{de_1pt}The following differential equations hold:
\ben
z\d_i \f_\tau^a = (\d_i\bullet)\f_\tau^a,\quad 
\d_\gl\f_\tau^{a} = -z^{-1}\f_\tau^a ,\quad
(z\d_z+\gl\d_\gl + E)\f_\tau^a = (\mu-1/2)\f_\tau^a + \frac{1}{k}\f_\tau^\phi.
\een
\end{lemma}
\proof
We refer to these three equations as the first, the second, and the last.  

The second equation is equivalent to $\d_\gl I_a^{(n)} = I_a^{(n+1)}.$
For $n\geq 0$ this equation holds by definition. Assume that 
$n=-p<0.$ 
Look at the definition \eqref{def:1pt}. The derivative with respect to $\gl$
of 
\beq\label{pt:integrand}
\frac{1}{p!}\int_{[x_a]} d^{-1}\((\gl-f_\tau)^p \omega \) =: \int_{[x_a]}\omega_p
\eeq
is a sum of two terms:
$\int_{[x_a]} d\omega_p/df_\tau + \int_{[x_a]} \d_\gl \omega_p.$ The first
term vanishes because $d\omega_p$ contains a factor of $(\gl-f_\tau)^p$ and 
$f_\tau(x_a)=\gl.$ The second term is precisely $\int_{[x_a]} \omega_{p-1}$.
The second equation is proved. 

For the last equation we will prove that the homogeneity properties of the Frobenius
structure implies the equality between the coefficients in front of $z^{-n},$ $n\geq 0.$
Then the equality between the positive powers of $z$ follows easily from the second equation. 
We have 
\ben
(\gl\d_\gl+E)\int_{[x_a]}\Phi_{\tau,\gl}^{(n)} = n 
\int_{[x_a]}\Phi_{\tau,\gl}^{(n)}, \quad
(\gl\d_\gl+E) \log x_a =1/k,
\een
because $\Phi_{\tau,\gl}^{(n)}(x_a)$ and $x_a$ are homogeneous in $\gl$
and $\tau$ of degrees $n$ and $1/k$ respectively. Now, to prove that the coefficients in front 
of $z^{-n}$ are equal, we only need to use that $\tau_i,1\leq i\leq k $ 
has degree $i/k,$ and $\mu(d\tau_i)= (1/2-i/k)d\tau_i,$ and 
$\tau_{k+j},1\leq j\leq m$ has degree $1-j/m$ and 
$\mu(d\tau_{k+j}) = (-1/2+j/m)d\tau_{k+j}.$

The first equation is equivalent to 
$\d_i I_a^{(n)} = -(\d_i\bullet)\d_\gl I_a^{(n)}.$ It is enough to
show that the later equation holds for all $n=-p,p\geq 2$ because for other
$n$ we just need to recall the second equation. 
Recall the definition \eqref{def:1pt} of $I_a^{(-p)}.$ We need to 
prove that for each $i$ and $j$ the differential operators 
$\d_i\d_j + \sum_l A_{ij}^l\d_l\d_\gl$ ($A_{ij}^l$ are the structure 
constants of $\bullet$)  annihilate the function \eqref{pt:integrand}.  
Note that
\ben
\(\d_i\d_j +\sum_l A_{ij}^l\d_l\d_\gl\)\int_{[x_a]}\omega_p = \int_{[x_a]}
\(\d_i\d_j + \sum_l A_{ij}^l\d_p\d_\gl\)\omega_p.
\een
In general the above formula will have more terms. However, in our case,
the integrands of those terms contain factors of $(\gl-f_\tau)^{p}$ or
$(\gl-f_\tau)^{p-1}$ and so they vanish on the one point cycle $[x_a].$ 
On the other hand we have 
\ben
&&
\(\d_i\d_j + \sum_l A_{ij}^l\d_l\d_\gl\)\frac{1}{p!}(\gl-f_\tau)^p\omega \\
&& 
=\[ -\frac{(\gl-f_\tau)^{p-1}}{(p-1)!}\,\frac{\d^2f_\tau}{\d\tau_i\d\tau_j} + 
\frac{(\gl-f_\tau)^{p-2}}{(p-2)!} (\d_if_\tau)(\d_jf_\tau) - 
\sum_l A_{ij}^l \frac{(\gl-f_\tau)^{p-2}}{(p-2)!}\d_pf_\tau\]\omega.
\een
Now we use that \eqref{primitive form} holds. Thus  the RHS of the last 
equation can be transformed into
\ben
-d \(G_{ij}\frac{(\gl-f_\tau)^{p-1}}{(p-1)!}\).
\een  
Now it is easy to finish the proof:
\ben
\(\d_i\d_j +\sum_l A_{ij}^l\d_l\d_\gl\)\int_{[x_a]}\omega_p = 
-\int_{[x_a]} G_{ij}\frac{(\gl-f_\tau)^{p-1}}{(p-1)!} = 0.
\een
\qed

\begin{corollary}\label{periods_de} 
Let $\gb\in H_1(\C^*,V_{\tau_0,\gl_0};\QQ)$ be any cycle. Then
\ben
z\d_i \f_\tau^\gb = (\d_i\bullet)\f_\tau^\gb,\quad 
\d_\gl\f_\tau^{\gb} = -z^{-1}\f_\tau^\gb ,\quad
(z\d_z+\gl\d_\gl + E)\f_\tau^\gb = (\mu-1/2)\f_\tau^\gb.
\een
\end{corollary}
\proof
Let $(\tau,\gl)\in \(M\times\C\)\backslash \Delta$. Using the Stokes' formula 
it is easy to see that there are one point cycles $x_a(\tau,\gl)$ and 
$x_b(\tau,\gl)$ such that $\f_\tau^\gb = \f_\tau^a-\f_\tau^b$, where the paths 
specifying the values respectively of $\f_\tau^a$ and $\f_\tau^b$
are possibly different from the path specifying the value of $\f_\tau^\gb.$ 
\qed


\subsection{Period vectors near $\gl=\infty$}
\label{sec:classical limit}

According to \leref{de_1pt}, $\f_\tau^a$ and $S_\tau$ satisfy the same 
differential equations with respect to $\tau$. In this section 
we will prove that in a neighborhood of $\gl=\infty$ the 
series $\f_\tau^a$ can be expanded as follows: 
$\f_\tau^a(\gl)=S_\tau\, \f_\infty^a(\gl)$, where 
$\f_\infty^a=\sum_n I_a^{(n)}(\infty,\gl)(-z)^n$ is independent of $\tau$
and each mode is a series of the following type,
\beq\label{periods_expansion}
\( \sum_{i=0}^{N_1} A_i \gl^{i}\)\log \gl  + 
\sum_{i\geq N_2} B_{i,1} \gl^{-i/k} + B_{i,2}\gl^{-i/m}, 
\eeq
where $A_i, B_{i,1},$ and $B_{i,2}$ are certain vectors in $H.$ 

We show that each period vector $I_a^{(n)}(\tau,\gl)$ expands 
in a neighborhood of $\gl=\infty$ as a series of type 
\eqref{periods_expansion} with coefficients in 
$H[\tau_1,\ldots,\tau_N,Qe^{\tau_N}]$. Then, to obtain 
$I_a^{(n)}(\infty,\gl),$ we just need to let
$\tau_1=\ldots=\tau_N=Qe^{\tau_N}=0.$ 

By definition,
\ben
\(I_\phi^{(-1)}(\tau,\gl),\d_i\) = -\d_i\int_\phi (\gl-f_\tau)\omega = \d_i t_k \int_\phi\omega=
\delta_{i,k}. 
\een
For the second equality we used that only the free term in $\gl-f_\tau$ contributes
to $\int_\phi$ and for the last one note that $t_k=\tau_k.$ 
Thus $I_\phi^{(-1)} = d\tau_k.$
The rest of the periods are uniquely determined from the differential 
equations in \coref{periods_de}.
Indeed, the second equation implies that $\d_\gl I_\phi^{(n)} =I_\phi^{(n+1)},$
thus $I_\phi^{(n)}=0$ for $n\geq 0$. Using the first equation we can express 
the differentiations $\d_i$ in terms of multiplication operators 
$\d_i\bullet.$ Thus the last equation assumes the following form 
\beq\label{recursion}
(\gl - E\bullet )I_\phi^{(-n+1)} = (\mu+n-1/2)I_\phi^{(-n)}.  
\eeq
In the basis $\d_1,\ldots,\d_N$ the Hodge grading operator $\mu$
is diagonal with entries 
\ben
\frac{1}{k}-\frac{1}{2},\ldots,  \frac{k-1}{k}-\frac{1}{2},
\frac{1}{2}-\frac{0}{m},\ldots, \frac{1}{2}-\frac{m}{m}.
\een
Thus $\mu+n-1/2$ is invertible for $n>1,$ i.e., $I_\phi^{(-n)}$ is 
uniquely determined from $I_\phi^{(-n+1)}.$ Moreover, it is clear 
that all modes $I_\phi^{(n)}$ will depend polynomially on $\gl$, $\tau$
and $Qe^{\tau_N}.$  

Similarly, the differential equations in \leref{de_1pt} imply that
the one-point periods satisfy the following recursive relation,
\beq\label{recursion_a}
(\gl - E\bullet )I_a^{(-n+1)} = (\mu+n-1/2)I_a^{(-n)} + 
\frac{1}{k}I_\phi^{(-n)}.
\eeq
Thus it suffices to prove that $I_a^{(-1)}(\tau,\gl)$ has an expansion 
of type \eqref{periods_expansion}. This is obvious because 
$I_a^{(-1)}(\tau,\gl)$ are obtained by integrating the 0-form
\beq\label{form}
\omega_1 = -\frac{x^k}{k} - \sum_{j=1}^{k-1}t_{k-j} \frac{x^{j}}{j}
+(\gl-t_k)\log x
+\sum_{j=1}^{m-1} \frac{t_{k+j}}{j} \(\frac{Qe^{t_N}}{x}\)^j 
+\frac{1}{m}\(\frac{Qe^{t_N}}{x}\)^m,
\eeq
over one-point cycles $x_a(\tau,\gl)$ -- the solutions of $f_\tau(x)=\gl$. Near
$\gl=\infty,$ the one-point cycle $x_a$ expands as a series in $\gl^{-1/k}$ or  
$\gl^{-1/m},$ depending on whether $1\leq a\leq k$ or $k+1\leq a\leq k+m$, 
whose coefficients depend polynomially on $\tau$ and $Qe^{\tau_N}.$ 

Now our goal is to compute $I_a^{(n)}(\infty,\gl)$ explicitly.  
Assume that $\gl$ is sufficiently 
large and that the path specifying the corresponding branch of 
$I_a^{(-1)}(\tau,\gl)$ is such that the one-point cycles (i.e. the 
solutions to $f_\tau(x)=\gl$) split into two groups: $x_1,\ldots,x_{k}$ close to 
$\infty$ and $x_{k+1},\ldots,x_{k+m}$ close to $0,$ and so the expansions
of $\log x_i$ coincide with the ones given in Section 
\ref{sec:flat_structure}. 
\begin{lemma}\label{periods_limit}
The vector-valued functions $\f_\infty^\phi,$ $\f_\infty^a,$ and $\f_\infty^b$
are given respectively by formulas \eqref{vector:phi}, \eqref{vector:a},
and \eqref{vector:b}. 
\end{lemma} 
\proof
We prove that $\f_\infty^b$ is given by \eqref{vector:b}. The other two 
cases are similar.

The series $\f_\infty^b$ is obtained from $\f_\tau^b$ by letting 
$\tau_1=\ldots=\tau_N=Qe^{\tau_N}=0.$ As explained
above, the equations $-z\d_\gl \f_\tau^b = \f_\tau^b$ and  
\eqref{recursion_a} determine $\f_\tau^b$ uniquely from $I_b^{(-1)}.$ 
Thus the same is true for $\f_\infty^b,$ i.e., it is uniquely 
determined from $I_b^{(-1)}(\infty,\gl)$ and the differential equations
(the second of the following two equations is obtained from 
\eqref{recursion_a} by letting $\tau_1=\ldots=\tau_N=Qe^{\tau_N}=0$):
\ben
(-z\d_\gl)\f_\infty^b =\f_\infty^b,\quad 
(z\d_z +\( \gl-\rho\)\d_\gl)\f_\infty^b =
(\mu -1/2)\f_\infty^b  + \frac{1}{k}\f_\infty^\phi ,    
\een
where $\rho$ is the classical multiplication by $(1/k+1/m)m\d_N.$ 
It is straightforward 
to check that the RHS of \eqref{vector:b} satisfies the
above differential equations. Thus it remains to prove that
\beq\label{I}
I_b^{(-1)}(\infty,\gl) = -\log (\gl^{1/m}Q^{-1}) d\tau_k - 
\sum_{j=1}^{m} 
\frac{\gl^{j/m}}{j/m} \d_{k+m-j}.
\eeq
By definition, 
\beq\label{xb}
(I_b^{(-1)}(\tau,\gl),\d_i)=-\d_i\int_{[x_b]} \omega_1 = 
\int_{[x_b]} \d_iF\frac{d\omega_1}{dF} - \int_{[x_b]}\d_i\omega_1.
\eeq
The first integral vanishes because the integrand contains a factor of
$(\gl-f_\tau)$. The
second one can be expanded into a series in $\log \gl$ and $\gl^{-1/m}.$
To obtain $I_b^{(-1)}(\infty,\gl)$ we need to ignore
all terms which depend on $\tau.$  Note that
$x_b = Qe^{t_N} \gl^{-1/m}+$(terms that depend on $\tau$). Thus our 
integral will have terms independent of $\tau$ only if 
$i=k+j,$ $0\leq j\leq m.$ When $i=k$ we get
\ben
-\int_{[x_b]}\d_k\omega_1=\log x_b = \log (Q\gl^{-1/m}) + \ldots,
\een
where the dots indicate terms depending on $\tau.$ Thus we obtain the 
logarithmic term in \eqref{I}. Furthermore, let $i=k+j,\ 1\leq j\leq m-1.$ 
Then, using the change of the variables \eqref{flat_coord} and \eqref{coord_flat}, 
we get (recall also \eqref{form})
\ben
-\int_{[x_b]}\d_{k+j}\omega_1  = 
\sum_{p=1}^{m-1} \frac{1}{p}(-\d_{k+j} t_{k+p}) 
\( \frac{Qe^{t_N}} {x_b} \)^p +\ldots = -\frac{1}{j} 
\gl^{j/m} +
\ldots,  
\een 
where the dots stand for terms which depend on $\tau.$ 
Finally, let $i=k+m=N.$ Recall that $\tau_N=mt_N.$ Thus we have 
\ben
-\int_{[x_b]}\d_{N}\omega_1  = -\int_{[x_b]}\frac{1}{m}\d_{t_N}\omega_1  = 
\frac{1}{m}
\(\frac{Qe^{t_N}}{x_b}\)^m+\ldots = -\frac{1}{m}\gl +\ldots . 
\een
The lemma is proved.
\qed


\section{Phase factors}\label{phasefactor}

The proof of Theorems \ref{t1} and \ref{t2} amounts to conjugating the vertex
operators $\Gamma_\infty^a$ and $\Gamma_\tau^a$ respectively with 
the symplectic transformations $S$ and $\Psi R.$ 
 
According to \cite{G2}, formula (17), conjugation of vertex operators
by lower-triangular transformations $S\in \L^{(2)}\lieGL(H)$ is given 
by the following formula:
\beq
\label{conjugation:s}
\widehat S e^{\hat \f}\hat S^{-1}  = 
e^{W(   \f_+  ,  \f_+   )/2}
e^{      (S \f) \sphat                      },
\eeq 
where $+$ means truncating the terms corresponding to the 
negative powers of $z$ and the quadratic form 
$W(\q,\q)=\sum (W_{kl}q_l,q_k)$ is defined by
\beq
\label{conjugation:w}
\sum_{k,l\geq 0}W_{kl}w^{-k}z^{-l} = \frac{S^*(w)S(z)-1}{w^{-1}+z^{-1}}.
\eeq
Similarly, for upper-triangular 
transformations $R\in \L^{(2)}\lieGL(H)$ we have (\cite{G2}, section 7):
\beq\label{conjugation:R}
\widehat R^{-1} e^{\hat \f} \widehat R = 
e^{V \f_-^2/2}\(e^{R^{-1}\f}\)\sphat,
\eeq 
where $-$ means truncating the non-negative powers of $z,$ $$\f_-=
\sum_{k\geq 0} (-1)^{-1-k}(f_{-1-k},q_k)$$ is
interpreted (via the symplectic form) as a linear function in $\q,$ and 
$V(\d,\d) = \sum (V_{kl}\phi^a,\phi^b)\d_{q_{l,a}}\d_{q_{k,b}}$ is 
a second order differential operator whose coefficients are defined by
\beq\label{conjugation:v}
\sum_{k,l\geq 0} V_{kl}w^kz^l = \frac{1-R(w)R^*(z)}{w+z}.
\eeq
In this section we will 
compute the {\em phase factors} $W(\f_+,\f_+)$ and $V\f_-^2$ respectively 
for the vertex operators $\Gamma_\infty^a$ and $\Gamma_\tau^\gb,$ 
where $\gb$ is a Lefschetz thimble.


\subsection{The phase form}
 
The phase factors can be conveniently expressed in terms of 
the so-called {\em phase form} $\W_{a,b}$ -- a multi-valued 
meromorphic section of $\pi^*(T^*M)$ (here $\pi:M\times\C\rightarrow M$ 
is the projection) with poles along the discriminant $\Delta$, 
defined by
\ben
\W_{a,b}(\tau,\gl) = \sum_{i=1}^N 
(I_a^{(0)}(\tau,\gl),\d_i\bullet I_b^{(0)}(\tau,\gl) )d\tau_i = 
I_a^{(0)}(\tau,\gl)\bullet I_b^{(0)}(\tau,\gl).
\een 
In the last equality vectors and {\em co}vectors are identified 
via the residue pairing.

\begin{lemma}\label{phase_form:primitive} Let $\tau\in M$ be a semi-simple 
point, and $x_a, x_b\in V_{\tau,\gl}$ two one-point cycles.  
Then
\ben 
\W_{a,b}(\tau,\gl) = 
\begin{cases} 
d^M \log\, (x_a-x_b) & \mbox{ if } a\neq b, \\ 
-d^M \log\,\( f_\tau'(x_a)\) & \mbox{ if } a=b,
\end{cases}
\een 
where $d^M$ is the de Rham differential on $M$ and $'$ means derivative 
with respect to $x.$
\end{lemma}
\proof
The idea is to express both sides in terms of canonical coordinates. 
Let $q_1,\ldots,q_N$ be the critical points of $f_\tau$ and $u_i$ are the
corresponding critical values, i.e.,
\ben
x^{m+1}f_\tau'(x) = k(x-q_1)\ldots(x-q_N),\quad u_i = f_\tau(q_i).
\een
That $\tau$ is a semi-simple point means that the critical values $u_i$
form a coordinate system on $M$ in a neighborhood of $\tau.$ Moreover, in 
that coordinate system the product and the metric assume the following forms:
\ben
\d/\d u_i\bullet\d/\d u_j = \delta_{ij}\d/\d u_i,\quad 
(\d/\d u_i,\d/\d u_i)={\delta_{ij}}/{\Delta_i},
\een 
where $\Delta_i$ is the Hessian of $f_\tau$ with respect to the volume form 
$\omega$ at the critical point $q_i,$ i.e., $\Delta_i =(\d_y^2 f_\tau)(q_i)$
where $y:=\log x$ is a unimodular coordinate on a neigborhood of $q_i$ in $\C^*.$   
Furthermore, introduce an auxiliary function 
\ben
\Phi:M\times \C^*\rightarrow \C,\quad \Phi(\tau,x) = f_\tau(x).
\een
Then the period vector $I_a^{(0)}$ can be written as follows
\beq\label{period:can}
I_a^{(0)} = -d^M\int_{[x_a]} d^{-1}\omega  = 
\sum_{i=1}^N \int_{[x_a]}\d_{u_i}\Phi\frac{\omega}{d\Phi} du_i = 
\sum_{i=1}^N\int_{[x_a]} \frac{\d_{u_i}\Phi}{\d_y\Phi}du_i,
\eeq
where $d$ is the de Rham differential on $\C^*.$

Borrowing an argument from the proof of Lemma 4.5 in
\cite{D}, we will express the partial derivative $\d_{u_i}\Phi$ in terms
of partial derivatives of $\Phi$ with respect to $y.$
By definition  $u_i(\tau) = \Phi(\tau,q_i).$ Applying $\d/\d u_j$
to both sides and using the chain rule, we get
\ben
\delta_{ij} = \d_{u_j}\Phi(q_i) + (\d_x\Phi)(q_i)\d_{u_j}q_i = 
\d_{u_j}\Phi(q_i)+f_\tau'(q_i)\d_{u_j}q_i = \d_{u_j}\Phi(q_i),
\een 
where we suppressed the dependence of $\Phi$ on $\tau.$
On the other hand $x^m\d_{u_j}\Phi$ is a polynomial in $x$ of degree
$N-1$ and the above formula shows that the value of this polynomial
at $q_i,$ $1\leq i\leq N$ is $q_i^m \delta_{ij}.$ Therefore the  Lagrange 
interpolation formula yields
\ben
\d_{u_j}\Phi(\tau,x) =\frac{q_j}{x-q_j}\,
\frac{\d_y\Phi(\tau,x)}{\d_y^2\Phi(\tau,q_j)}=
\frac{q_j}{\Delta_j(x-q_j)}\,\d_y\Phi(\tau,x). 
\een  
Hence formula \eqref{period:can} transforms into
\ben
I_a^{(0)}(\tau,\gl) = \sum_{i=1}^N \frac{q_i}{\Delta_i(x_a-q_i)}\, du_i.
\een
Using that $du_i\bullet du_j=\delta_{ij}\Delta_idu_j$ we find that in 
canonical coordinates the phase form is given by the following 
formula
\ben
\W_{a,b} = \sum_{i=1}^N \frac{q_i^2}{\Delta_i(x_a-q_i)(x_b-q_i)}\, du_i.
\een
Assume that $x_a\neq x_b.$ Then 
\beq\label{primitive}
d^M\log(x_a-x_b) = \sum_{i=1}^N\frac{\d_{u_i}x_a-\d_{u_i}x_b}{x_a-x_b}du_i. 
\eeq 
By definition $\gl = \Phi(\tau,x_a).$  Apply $\d/\d u_i$ 
to both sides and solve the resulting identity with respect to $\d_{u_i}x_a$, we find:
\ben
\d_{u_i}x_a=-\frac{\d_{u_i}\Phi(x_a)}{\d_x\Phi(x_a)} =
-\frac{x_a\d_{u_i}\Phi(x_a)}{\d_y\Phi(x_a)} = 
-\frac{x_a q_i}{\Delta_i(x_a-q_i)}. 
\een
Substitute this in \eqref{primitive}. After a short 
simplification we get precisely the formula for $\W_{a,b}.$

The case when $a=b$ is similar. We have
\ben
-d^M\log (\d_x\Phi(x_a)) = 
-\sum_{i=1}^N\frac{1}{\d_x\Phi(x_a)}
\[(\d_x\d_{u_i}\Phi)(x_a) + (\d_x^2\Phi)(x_a)\d_{u_i}x_a\]du_i.
\een 
The first term in the brackets equals to
\ben
\left( \d_x \frac{q_i}{\Delta_i(x-q_i)}\d_y\Phi \right)(x_a) = 
-\frac{q_ix_a}{\Delta_i(x_a-q_i)^2}\d_x\Phi(x_a) + 
\frac{q_i}{\Delta_i(x_a-q_i)}(\d_x\d_y\Phi)(x_a).
\een
The second one equals to
\ben
-(\d_x^2\Phi)(x_a)\frac{q_i x_a}{\Delta_i(x_a-q_i)}.
\een
Using that $(\d_x\d_y\Phi)(x_a) = \d_x\Phi(x_a) + x_a\d_x^2\Phi(x_a)$
we find that the sum of the two terms is 
\ben
-\frac{q_ix_a}{\Delta_i(x_a-q_i)^2}\d_x\Phi(x_a)+ 
\frac{q_i}{\Delta_i(x_a-q_i)}\d_x\Phi(x_a) = 
-\frac{q_i^2}{\Delta_i(x_a-q_i)^2}\d_x\Phi(x_a).
\een
Therefore $d^M\log (\d_x\Phi(x_a))=\W_{a,a}.$
\qed

\begin{corollary}\label{co:phase_form}
Assume the notations from \leref{phase_form:primitive}. Then
\ben
\(I_a^{(0)}(\tau,\gl),I_b^{(0)}(\tau,\gl)\) = 
\begin{cases}
-\d_\gl \log\, (x_a-x_b) & \mbox{ if } a\neq b, \\ 
\d_\gl \log\,\( f_\tau'(x_a)\) & \mbox{ if } a=b.
\end{cases}
\een  
\end{corollary}
\proof
It follows from the definitions that the one-point cycles $x_a(\tau,\gl)$
depend only on the difference $t_k-\gl,$ i.e., $\d_\gl x_a = -\d_k x_a.$
Recall \leref{phase_form:primitive} and pair the differential 1-forms 
with the vector field $\d_k.$
\qed


\subsection{Phase factors near $\infty$}


Assume that $\gl$ is sufficiently large and that the path specifying
the branches of the period vectors $I_i^{(0)}$ is such that 
$\log x_i$ have the same expansions as in Section \ref{sec:flat_structure}.
For each $i,$ $1\leq i\leq N$ put $\f^i = \f_\infty^i$ to avoid cumbersome
notations. Using  definition \eqref{conjugation:w} of the quadratic form 
$W_\tau$ and that $\d_\gl I_a^{(n)} = I_a^{(n+1)}$ we get 
{\allowdisplaybreaks
\ben
&&
\frac{d}{d\xi} W_\tau(\f^i_+(\xi),\f^i_+(\xi)) = 
-\sum_{k,l\geq 0} 
\( (W_{k,l-1}+W_{k-1,l})(-1)^l I_i^{(l)}(\xi), 
                        (-1)^k I_i^{(k)}(\xi) ) \) \\
& = &
-\sum_{k,l\geq 0} 
\( S_l(-1)^l I_i^{(l)}(\xi), 
                        S_k(-1)^k I_i^{(k)}(\xi) ) \) 
+\( I_i^{(0)}(\xi), I_i^{(0)}(\xi) \)  \\
& = &
-\(I_i^{(0)}(\tau,\xi),I_i^{(0)}(\tau,\xi)\) +
\( I_i^{(0)}(\xi), I_i^{(0)}(\xi)\) .
\een}
There are two cases: $i=a,$ $1\leq a\leq k$ or $i=b,$ $k+1\leq b\leq k+m$
which correspond respectively  to $x_i$ being close to $x=\infty$ or 
$x=0.$ In the first case $\f^a_+=0$ at $\xi=\infty.$ Thus we get
\ben
W_\tau(\f^a_+(\gl),\f^a_+(\gl)) = \int^\infty_\gl \[ 
(I_a^{(0)}(\tau,\xi),I_a^{(0)}(\tau,\xi)) - 
( I_a^{(0)}(\infty,\xi), I_a^{(0)}(\infty,\xi)) \]d\xi,
\een 
where the integration path is a ray starting at $\gl$ and 
approaching $\gl=\infty.$
According to \leref{periods_limit}, 
\ben 
I_a^{(0)}(\infty,\xi)= \sum_{i=1}^{k-1} \xi^{i/k-1}\d_i + 
\frac{1}{k}\xi^{-1}d\tau_k.
\een
Thus we get the following formula:
\beq\label{phase_factor:a}
W_\tau(\f^a_+,\f^a_+) = \int^\infty_\gl \[
(I_a^{(0)}(\tau,\xi),I_a^{(0)}(\tau,\xi)) - \frac{k-1}{k} \xi^{-1} \]d\xi.
\eeq
In the second case $\f_+^b=-\d_k$ at $\xi=\infty,$ therefore 
the phase factor $W_\tau(\f_+^b(\gl),\f_+^b(\gl))$ is given by
the following integral:
\ben 
W_\tau(\d_k,\d_k) +
\int^\infty_\gl \[
(I_b^{(0)}(\tau,\xi),I_b^{(0)}(\tau,\xi)) - 
( I_b^{(0)}(\xi,\infty), I_b^{(0)}(\xi,\infty)) \]d\xi,
\een 
where the integration path is the same as above. 
According to \leref{periods_limit},
\ben
I_b^{(0)}(\infty,\xi) = -\sum_{j=1}^{m-1}\xi^{j/m-1}\d_{k+m-j} -\d_k-
\frac{1}{m} \xi^{-1}d\tau_k. 
\een 
The term $W_\tau(\d_k,\d_k)=(S_1\d_k,\d_k)$ can be computed as follows:
\ben
\d_i (S_1\d_k,\d_k) = (\d_iS_1\d_k,\d_k) = (\d_i\bullet\d_k,\d_k)=
(\d_k,\d_i) = \frac{1}{m}\delta_{i,N}.
\een
On the other hand $S_1=0$ when $\tau=0$, so $(S_1\d_k,\d_k)=\tau_N/m = t_N.$
Thus 
\beq\label{phase_factors:b}
W_\tau(\f^b_+,\f^b_+) = t_N + \int^\infty_\gl \[
(I_b^{(0)}(\xi,\tau),I_b^{(0)}(\xi,\tau)) - 
\frac{m+1}{m}\xi^{-1} \]d\xi.
\eeq


\subsection{Phase factors near a critical point}

Let $\tau\in M$ be such that $f_\tau$ is a Morse function. Assume that $\gl$ 
is close to a critical value $u_i,$ $1\leq i\leq N,$ of $f_\tau$ and that
$\gb\in H_1(\C^*,V_{\tau,\gl};\Z)$ is a Lefschetz thimble corresponding to the path 
$C_i,$ where $C_i$ is the straight segment from $\gl$ to $u_i.$ 
Let $\d\gb=[x_a]-[x_b],$ where the one-point cycles are represented by some 
$x_a,x_b\in V_{\tau,\gl}.$ 
  
\begin{lemma}\label{phase_factor}
Let $V$ be the quadratic form \eqref{conjugation:v}. Then
\beq\label{phase}
V\f_-^2 = 
-\lim_{\ge\rightarrow 0}
\int_{\gl}^{u_i+\ge}
\(
\( I_{\gb/2}^{(0)}(\tau,\xi) , I_{\gb/2}^{(0)}(\tau,\xi)\)
-\frac{1}{ 2(\xi-u_i) } 
\) d\xi , 
\eeq
where $\f:=\f_\tau^{\gb/2}$ and the limit is taken along the 
path $C_i.$ 
\end{lemma}
\proof 
The proof is taken from \cite{G2}, page 490. 
When 
$\f=\sum_{k\in\Z} I_{\gb/2}^{(k)}(-z)^k$, we have 
$\f_-=\sum_{k\geq 0} \sum_{a=1}^N( I_{\gb/2}^{(-1-k)},\phi_a)q_{k}^a.$
Using $\d_\gl I_{\gb}^{(-1-k)}=I_{\gb}^{(-k)},$ we find  
\ben
\d_\gl V\f_-^2 
&=& \frac{1}{4}\sum_{k,l\geq 0} \d_\gl 
\(V_{k,l}I_{\gb}^{(-1-l)},I_{\gb}^{(-1-k)}\) = 
\frac{1}{4}\sum_{k,l\geq 0} 
\( [V_{k-1,l} + V_{k,l-1} ] I_{\gb}^{(-l)},I_{\gb}^{(-k)} \)\\
&=&
\frac{1}{4}(I_{\gb}^{(0)},I_{\gb}^{(0)})- 
\frac{1}{4}
(\sum_{l\geq 0}R_l^*I_{\gb}^{(-l)},\sum_{k\geq 0}R_k^*I_{\gb}^{(-k)})  .
\een
On the other hand, $R^*(z) = R^{-1}(-z)$ because $R\in \L^{(2)}\lieGL(H).$
Thanks to \leref{vector}, 
$\sum_{k\geq 0}R_k^*I_{\gb}^{(-k)}=
2{\bf 1}_i/\sqrt{2(\gl-u_i)}.$ Also $V\f_-^2 = 0$ at $\gl=u_i$ because 
$I^{(-1-k)}_{\gb_i} = 2[2(\gl-u_i)]^{k+1/2}({\bf 1}_i +\ldots )$ vanish
at $\gl=u_i.$ The lemma follows. 
\qed

\begin{lemma}\label{vop:invariant}
Let $\ga=s([x_a]+[x_b]),\ s\in \QQ.$ Then the period vectors 
$I_{\ga}^{(n)} := s( I_{a}^{(n)}+I_{b}^{(n)})$ are holomorphic in 
a neighborhood of $\gl=u_i$ for each $n\in \Z.$
\end{lemma}
\proof
The function $f_\tau$ is Morse in a neighborhood of $\gl=u_i.$ Let $y$ be 
a Morse coordinate, i.e., 
\beq\label{morse_coord}
f_\tau = u_i + \frac{y^2}{2}, 
\quad x=q_i(\tau) + a_1(\tau) y + a_2(\tau)y^2+\ldots , 
\eeq
where $q_i$ is a critical point of $f_\tau$ corresponding to the critical value 
$u_i.$ 
Therefore, in the Morse coordinate, the one-point cycles $x_a(\tau,\gl)$ and 
$x_b(\tau,\gl)$ are given by $y_{a/b} = \pm \sqrt{2(\gl-u_i)}.$ The last formula shows that the 
period vectors $I_{\ga}^{(n)}$ (which by definition are integrals of certain 0-forms
over the one-point cycles $[x_a]+[x_b]$) are single--valued near $\gl=u_i.$ Moreover, 
they are obtained from  $I_{\ga}^{(0)}$ by differentiating and antidifferentiating 
with respect to $\gl$. Thus it is enough to show that  $I_{\ga}^{(0)}$ is 
holomorphic near $\gl=u_i.$ 

On the other hand $I_{\ga}^{(0)} = -s \, d^M \, (\log x_a+ \log x_b ).$
The last expression is holomorphic near $\gl=u_i$ because it is singlevalued
(the analytical continuation around $u_i$ transforms $x_{a/b}$ into $x_{b/a}$) and 
$x_a,x_b\rightarrow q_i$ when $\gl\rightarrow u_i$ and $q_i\neq 0.$ 
\qed

\begin{lemma}\label{vop:composition}
Let $\ga = s([x_a]+[x_b]).$ Then the following formula holds:
\ben
\Gamma_\tau^{\ga + r\gb/2} = \exp\(r\int^\gl_u 
\(I_{\ga}^{(0)}(\tau,\xi),I_{\gb/2}^{(0)}(\tau,\xi)\) d\xi \) 
\Gamma_\tau^\ga \Gamma_\tau^{r\gb/2},
\een
where $r\in \QQ$ is an arbitrary number. 
\end{lemma}
\proof
See \cite{G2}, Proposition 4. 
\qed


\section{From ancestors to KdV}\label{tokdv}

Let $\tau\in M$ be an arbitrary semi-simple point. We want to
show that the ancestor potential $\A_\tau^M$ satisfies the 
HQE \eqref{HQE:ancestors}. A priori the vertex operators in 
\eqref{HQE:ancestors} depend on the choice of a path $C$ in 
$\C\backslash \{u_1,\ldots,u_N\}$ from $\gl_0$ to $\gl.$ 
We begin by showing that the HQE are independent of $C$.  

\subsection{Monodromy and vertex operators}

Assume that $C'$ is another path from $\gl_0$ to $\gl$, denote
the corresponding one-point cycles by $x_a'(\tau,\gl)$  and the 
associated vertex operator by $\Gamma_\tau^{a'}.$ Since $V$ is a
covering we have that $x_a' = x_{\sigma(a)}$ for some permutation 
$\sigma$ of $1,2,\ldots, N.$   Using the definition of 
the period vectors we find
\ben
\(I_{a'}^{(-n)}-I_{\sigma(a)}^{(-n)},\d_i\) = 
-\d_i \int_{[x_a']-[x_{\sigma(a)}]}d^{-1}\(\frac{1}{n!}(\gl- f)^n\omega\) = 
-\d_i \int_{r\phi} \frac{1}{n!}(\gl- f)^n\omega, 
\een
for some $r\in 2\pi i\Z.$ Thus 
$\f_\tau^{a'}=\f_\tau^{\sigma(a)} +r\f_\tau^\phi.$ Using that 
$I_\phi^{(n)}(\tau,\gl)=0$ for $n\geq 0$ we get 
$\Gamma_\tau^{a'}= \Gamma_\tau^{r\phi}\Gamma_\tau^{\sigma(a)}.$ It is
obvious that $c_{a'}=c_{\sigma(a)},$ therefore 
changing the path $C$ transforms the HQE \eqref{HQE:ancestors} into 
\ben
\(\Gamma_\tau^{\gd\#}\tensor \Gamma_\tau^\gd \) 
\(
\sum_{i=1}^N c_\tau^{\sigma(i)}
\(\Gamma_\tau^{r_i\phi}\tensor \Gamma_\tau^{-r_i\phi}\)
\(\Gamma_\tau^{\sigma(i)}\tensor \Gamma_\tau^{-\sigma(i)}\) 
\)
\(\T\tensor \T\)\  {d\gl}.
\een
The later expression, when computed at $\q',\q''$ such that 
$\hat w_\tau'-\hat w_\tau''\in \Z$, coincides with \eqref{HQE:ancestors}
thanks to the following Lemma.
\begin{lemma}{\label{OPE1}}
Let $r\in \C.$
The following identity between operators acting on $\B_H\tensor_\A\B_H$
holds:
\ben
\(\Gamma_\tau^{\delta \#}\tensor\Gamma_\tau^\delta\)
\(\Gamma_\tau^{r\,\phi}\tensor  \Gamma_\tau^{-r\,\phi}\) =
e^{(\hat w_\tau\tensor 1 - 1\tensor \hat w_\tau)\,r}
\(\Gamma_\tau^{\delta \#}\tensor\Gamma_\tau^\delta\).
\een
\end{lemma}
\proof
By definition   
\ben
\f_\tau^\phi(\gl) =  d\tau_k (-z)^{-1} + 
\sum_{n\geq 1} I_\phi^{(-1-n)}(\tau,\gl)(-z)^{-k-1}.
\een
Since $w_\infty =-d\tau_k z^{-1}$ and  $w_\tau = S_\tau w_\infty$, we get
$\f_\tau^{\phi}(\gl)\in w_\tau+z^{-1}\H_-.$ Thus
\ben
\Omega(v_\tau,\f_\tau^{\phi}(\gl)) = \Omega(v_\tau,w_\tau)=-1,\ \ 
\Omega(w_\tau,\f_\tau^{\phi}(\gl)) =0 \ .
\een

For all $f,g\in \H$ we have: 
\ben
e^{\hat f}e^{\hat g} = e^{[\hat f,\hat g]}e^{\hat g}e^{\hat f} = 
                       e^{\Omega(\hat f,\hat g)}e^{\hat g}e^{\hat f},
\een
because for linear Hamiltonians the quantization is a representation of
Lie algebras. 
Thus
{\allowdisplaybreaks
\ben
&&
\Gamma_\tau^\delta\Gamma_\tau^{r\,\phi}\\
& = &
\exp\left( 
(\widehat \f_\tau^{\phi}-\widehat w_\tau)(\ge\d_x) \right)
\exp \( \frac{x}{\ge}\widehat v_\tau\)\exp\(r\widehat\f_\tau^\phi\) \\
& = &
\exp\(r\widehat\f_\tau^\phi\)\exp\left( 
(\widehat \f_\tau^{\phi}-\widehat w_\tau)(\ge\d_x) \right)
\exp\((r x/\ge)\Omega(v_\tau,\f_\tau^\phi)\)\exp \( \frac{x}{\ge}\widehat v_\tau\)\\
& = &
e^{\hat w_\tau r}e^{ - rx/\ge} \Gamma_\tau^\delta
\een}
Similarly,
$
\Gamma_\tau^{\delta\#}\Gamma_\tau^{-r\,\phi} = 
e^{-\hat w_\tau\,r} \Gamma_\tau^{\delta\#} e^{x r/\ge}.
$
The Lemma follows. 
\qed

\subsection{Tame asymptotical functions}
The total ancestor potential has some special property
which makes the expression \eqref{HQE:ancestors} a formal series 
with coefficients meromorphic functions in $\gl.$ 

\begin{definition}[\cite{G2}]
{\em
An {\em asymptotical function} is by definition an expression 
\ben
\T = \exp \left(\sum_{g=0}^\infty \ge^{2g-2}\T^{(g)}(\t;Q)\right),
\een
where $\T^{(g)}$ are formal series in the sequence of vector 
variables $t_0,t_1,t_2,\ldots$  with coefficients in the Novikov
ring $\C[[Q]].$ Furthermore, $\T$ is called {\em tame} if
\ben
\left.\frac{\d}{\d t_{k_1}^{a_1}} \ldots \frac{\d}{\d t_{k_r}^{a_r}}
\right|_{\t=0} 
\T^{(g)} = 0 \quad \mbox{whenever} \quad
k_1+k_2+\ldots +k_r > 3g-3+r,
\een
where $t_k^a$ are the coordinates of $t_k$ with respect to 
$\{\phi_a\}, 1\leq a\leq N.$}
\end{definition}
According to \cite{G2}, Proposition 5, the total ancestor potential 
$\A_\tau^M$ is a tame asymptotical function.

Let $\T$ be a tame asymptotical function.  
The dilaton shift $\t(z)=\q(z)+z$ identifies $\T$ with an asymptotical element
of the Fock space $B_H$. Let 
$\Gamma_i^{\pm} = \exp\(\pm\sum I_i^{(n)}(\gl)(-z)^n\)\sphat$ be a finite
set of vertex operators, where $I_i^{(n)}$ are meromorphic functions. 
Consider the expression 
\beq\label{tameness}
\sum_i c^i(\gl)\(\Gamma_i^+\tensor\Gamma_i^-\)\(\T(\q')\otimes\T(\q'')\)d\gl, 
\eeq
where $c^i$ are meromorphic functions. 
According to \cite{G2}, Proposition 6, the tameness of $\T$ implies that
\eqref{tameness}, after the substitution $\q'=\x+\ge \y$, 
$\q''=\x-\ge\y$ and dividing by $\exp\(2\T^{(0)}(\x)/\ge^2\),$
expands into a power series in $\ge,$ $\x,$ and $\y$ whose 
coefficients depend polynomially on finitely many $I_i^{(n)}.$

In particular, \eqref{HQE:ancestors} can be interpreted as a formal series 
in $\ge$, $\x,$ and $\y$ (with $y_{0}^k$ excluded) with coefficients 
meromorphic functions in $\gl.$ The vertex operators could have 
poles only at the critical values $u_i, 1\leq i\leq N.$ Thus the 
regularity property of \eqref{HQE:ancestors} follows if we prove 
that there are no poles at $\gl=u_i.$   

\subsection { Proof of \thref{t2}}
Fix an arbitrary critical value $u_i.$ We need to show that 
\eqref{HQE:ancestors} does not have a pole at $\gl=u_i.$ Assume that
$\gl$ is sufficiently close to $u_i$ and let $C_i$ be the straight
segment from $\gl$ to $u_i.$ Then there are exactly two 
one-point cycles $x_a(\tau,\gl)$ and $x_b(\tau,\gl)$ which will coincide
after transported along $C_i$ towards the critical point. Note that
all period vectors $I_i^{(n)},$ $i\neq a$ or $b,$ are holomorphic
functions in a neighborhood of $\gl=u_i.$ Also, the vertex
operator $\Gamma_\tau^{\delta}$ is holomorphic (even polynomial) in
$\gl.$ Thus we need to show that 
\beq\label{HQE:ui}
 \(
\frac{1}{f_\tau'(x_a)}
\Gamma_\tau^{a}\tensor \Gamma_\tau^{-a} +
\frac{1}{f_\tau'(x_b)}
\Gamma_\tau^{b}\tensor \Gamma_\tau^{-b} 
\)
\(\A_\tau\tensor \A_\tau\)\  {d\gl}
\eeq
has no poles in a neighborhood of $u_i.$ 

The 1-point cycle $x_a$ can be splitted into 
$x_a= (x_a-x_b)/2 + (x_a+x_b)/2$. Using \leref{vop:composition}
we get
\ben
\Gamma_\tau^a\tensor\Gamma_\tau^{-a} = e^{K_a}
\(\Gamma_\tau^{\ga}\tensor\Gamma_\tau^{-\ga}\)  
\(\Gamma_\tau^{\gb/2}\tensor\Gamma_\tau^{-\gb/2}\), 
\een
where $\ga=(x_a+x_b)/2,$ $\gb$ is the Lefschetz thimble corresponding
to the critical point $u_i,$ and  
\ben
K_a = \int_{u_i}^\gl \(I_{\ga}^{(0)},I_{\gb}^{(0)}\)d\xi.
\een
Similarly,  
\ben
\Gamma_\tau^b\tensor\Gamma_\tau^{-b} = e^{K_b}
\(\Gamma_\tau^{\ga}\tensor\Gamma_\tau^{-\ga}\)  
\(\Gamma_\tau^{-\gb/2}\tensor\Gamma_\tau^{\gb/2}\), 
\een
where $K_b = -K_a.$

Furthermore, we recall formula \eqref{conjugation:R}: 
\ben
\(\Gamma_\tau^{\pm\gb/2}\tensor\Gamma_\tau^{\mp\gb/2}\)
\(\widehat\Psi \widehat R\tensor\widehat\Psi \widehat R\) = 
e^W \(\widehat \Psi \widehat R\tensor\widehat \Psi \widehat R\) 
\(\Gamma_{u_i}^{\pm}\tensor\Gamma_{u_i}^{\mp}\),
\een
where, according to \leref{phase_factor},
\ben
W=\lim_{\ge\rightarrow 0} \int_{u_i+\ge}^\gl \(
\(I_{\gb/2}^{(0)}(\tau,\xi),I_{\gb/2}^{(0)}(\tau,\xi)\)-\frac{1}{2(\xi-u_i)} \)
d\xi ,
\een
and according to \leref{vector}, 
\ben
\Gamma_{u_i}^{\pm} = 
\exp \( \pm\sum_{n\in\Z} (-z\d_\gl)^n \frac{e_i}{\sqrt{2(\gl-u_i)}}\).
\een
Finally, the periods $I_{\ga}^{(n)},\ n\in \Z,$ are holomorphic near
$\gl=u_i$ (see \leref{vop:invariant}).  Therefore,  to prove that 
\eqref{HQE:ui} is holomorphic near
$\gl=u_i,$ it is enough to show that the 1-form 
\beq \label{HQE:uii}
\(C_a \Gamma_{u_i}^{+}\tensor \Gamma_{u_i}^{-} + 
C_b \Gamma_{u_i}^{-}\tensor \Gamma_{u_i}^{+}\) 
\( \D_{\rm pt}\tensor \D_{\rm pt}\)d\gl ,
\eeq
is analytic in $\gl.$ Here 
\ben
C_{a/b} = \exp\(  K_{a/b}+W- \log {f_\tau'(x_{a/b})}\),
\een
i.e.,
\ben
\log C_{a/b} = - \log {f_\tau'(x_{a/b})} + \lim_{\ge \rightarrow 0} 
\int_{u_i+\ge}^\gl \( \pm (I_{\ga}^{(0)},I_\gb^{(0)}) + (I_{\gb/2}^{(0)},I_{\gb/2}^{(0)})
-\frac{1}{2(\xi-u_i)} \) d\xi ,
\een
where the periods in the integrand are computed at the point $(\tau,\xi).$ 
If we change $C_{a/b}$ by adding $(I_{\ga}^{(0)} ,I_\ga^{(0)})$ to the integrand of the RHS, 
then the expression \eqref{HQE:uii} will change by an
invertible holomorphic factor. Thus we can assume that the integrand is given by
{\allowdisplaybreaks
\ben
(I_{\ga}^{(0)} ,I_\ga^{(0)})
\pm (I_{\ga}^{(0)},I_\gb^{(0)}) + (I_{\gb/2}^{(0)},I_{\gb/2}^{(0)})
-\frac{1}{2(\xi-u_i)}  = 
(I_{a/b}^{(0)}, I_{a/b}^{(0)}) -\frac{1}{2(\xi-u_i)} .
\een}
Now recall \coref{co:phase_form} 
\ben
\log C_{a/b} & = & - \log {f_\tau'(x_{a/b})} + \lim_{\ge \rightarrow 0}
\int_{u_i+\ge}^\gl d \( \log \frac{f_\tau'(x_{a/b})}{ \sqrt{2(\xi-u_i)}}\) \\
& = &
-\log \sqrt{2(\gl-u_i)} + \lim_{\xi\rightarrow u_i} 
\log \frac{f_\tau'(x_{a/b})}{ \sqrt{2(\xi-u_i)}}, 
\een
where, in the last expression, the one-point cycles are computed at the point 
$(\tau,\xi).$  We will show that the above limit is 
$\log (\pm \sqrt{f_\tau''(q_i)}),$ where $q_i$ is a critical point of $f_\tau$ with 
critical value $u_i.$ Indeed, let us expand $f_\tau$ in a Taylor's series 
about $x=q_i$:
\ben
f_\tau(x) = u_i + \frac{1}{2!} f_\tau''(q_i) (x-q_i)^2 + \ldots .
\een
From this we find that $x_{a/b}(\tau,\xi)$ expands into a series in $\xi-u_i:$ 
\ben
x_{a/b} = q_i \pm \frac{1}{\sqrt{f_\tau''(q_i)}}\, \sqrt{2(\xi-u_i)} + \ldots,
\een
where the dots stand for higher order terms in $\xi-u_i.$ Moreover, differentiating
the Taylor's expansion in $x$ and then substituting $x=x_{a/b}$ yield
\ben
f_\tau'(x_{a/b}) = \pm \sqrt{f_\tau''(q_i)}\sqrt{2(\xi-u_i)} + \ldots .
\een
Note that $\log C_a - \log C_b = -\log (-1),$ i.e., $C_a = -C_b.$ 
Therefore the expression \eqref{HQE:uii}, up to an invertible holomorphic factor, equals to
\ben
\(
\frac{1}{\sqrt{2(\gl-u_i)}}\, \Gamma_{u_i}^{+} \tensor \Gamma_{u_i}^{-} - 
\frac{1}{\sqrt{2(\gl-u_i)}}\, \Gamma_{u_i}^{-} \tensor \Gamma_{u_i}^{+} \)
\(\D_{\rm pt} \tensor \D_{\rm pt} \).
\een 
This is precisely \eqref{HQE:kdv} and is therefore holomorphic at $\gl=u_i.$ 
\qed


\section{From ancestors to descendents}
\label{sec:From ancestors to descendents}


In this section we prove \thref{t1}. Recall that 
the descendent and the ancestor potentials are related by 
$\D^M = C_\tau \widehat S_\tau^{-1}\A_\tau,$ where $C_\tau$ is some 
constant and $S_\tau$ is the calibration of $M.$ The action of
$\widehat S_\tau^{-1}$ on the Fock space $B_H$ is given by
the following formula (\cite{G3}, Proposition 5.3):
\beq\label{action:s}
\( \widehat S_\tau^{-1} {\mathcal{G}}\) (\q) = 
e^{W_\tau(\q,\q)/2\ge^2} {\mathcal{G}} ([S_\tau \q]_+),
\eeq
where $W_\tau$ is a quadratic form defined by \eqref{conjugation:w} and 
$[\ ]_+$ means truncating the terms with negative powers of $z.$ 

On the other hand $S_\tau$ acts on the set of vertex operators  
$\Gamma_\infty^i,\ 1\leq i\leq N,$ and $\Gamma_\infty^\delta$ by 
conjugation. 
\begin{lemma}\label{conjugation:a} 
The following formula holds:
\ben
c_\infty^i \( \Gamma_\infty^i\tensor \Gamma_\infty^{-i} \) 
\(\widehat S_\tau^{-1}\tensor \widehat S_\tau^{-1}\) = 
\(\widehat S_\tau^{-1}\tensor \widehat S_\tau^{-1}\) c_\tau^i(\gl)
\( \Gamma_\tau^i\tensor \Gamma_\tau^{-i} \).
\een
\end{lemma}
\proof
We prove the case $1\leq i\leq k.$ The other case $k+1\leq i\leq N$ is similar.  

According to formulas \eqref{conjugation:s} and \eqref{phase_factor:a} we have
\ben
c_\infty^i \( \Gamma_\infty^i\tensor \Gamma_\infty^{-i} \) 
\(\widehat S_\tau^{-1}\tensor \widehat S_\tau^{-1}\) = 
C_i  
\(\widehat S_\tau^{-1}\tensor \widehat S_\tau^{-1}\) 
\( \Gamma_\tau^i\tensor \Gamma_\tau^{-i} \),
\een
where 
\ben
\log C_i= \frac{1-k}{k}\log \gl -\log k + \int_\gl^\infty 
\[\(I_i^{(0)}(\tau,\xi), I_i^{(0)}(\tau,\xi)\) - \frac{k-1}{k}\xi^{-1}\]d\xi.
\een
On the other hand, using \coref{co:phase_form}, we find that the above integral equals
\ben
\lim_{\xi\rightarrow\infty} \ \(\log  \frac{f_\tau'( x_i(\xi)) }{ \xi^{(k-1)/k} }\)
-\log f_\tau'(x_i) + \frac{k-1}{k}\log \gl.
\een 
Thus it remains to show that the above limit is $\log k.$
Indeed, near $\xi = \infty,$ we have $x_i(\xi) = \xi^{1/k} + ...,$ where here and 
further the dots stand for lower order terms. 
Hence $f_\tau'(x_i) = k x_i^{k-1} + \ldots = k \xi^{(k-1)/k} + \ldots .$ The lemma 
follows. 
\qed

\begin{lemma}\label{conjugation:d}The following formula holds:
\ben
\( \Gamma_\infty^{\delta \#}\tensor \Gamma_\infty^{\delta} \) 
\(\widehat S_\tau^{-1}\tensor \widehat S_\tau^{-1}\) = 
e^{t_N x^2/(2\ge)}\(\widehat S_\tau^{-1}\tensor \widehat S_\tau^{-1}\)
\( \Gamma_\tau^{\delta \#}\tensor \Gamma_\tau^{\delta} \)e^{t_N x^2/(2\ge)}.
\een
\end{lemma} 
\proof
Note that
$w_\tau=S_\tau \, w_\infty $ and $v_\tau = S_\tau v_\infty.$
Using \eqref{conjugation:s} we compute: 
\ben
&&
\widehat S_\tau^{-1} \Gamma_\tau^\gd \widehat S_\tau = 
\widehat S_\tau^{-1} \exp \(
\(\f_\tau^{\phi} -w_\tau\)\ge\d_x\)\sphat \ 
\exp\( x\hat v_\tau/\ge \) \sphat 
\hat S_\tau  \\
&&
=\exp\( 
(\f_\infty^{\phi} -w_\infty )\ge\d_x  \)\sphat
\exp\(
xv_\infty/\ge\)\sphat 
e^{-W_\tau(v_\infty,v_\infty)\frac{x^2}{2\ge^2}}.  
\een
On the other hand $W_\tau(v_\infty,v_\infty)=W_\tau(\d_k,\d_k) = t_N.$ Thus 
$\widehat S_\tau^{-1} \Gamma_\tau^\gd \widehat S_\tau=
\Gamma_\infty^\delta e^{-t_Nx^2/(2\ge^2)}.$ Similarly, 
$\widehat S_\tau^{-1} \Gamma_\tau^{\gd\#} \widehat S_\tau=
e^{-t_N x^2/(2\ge^2)}\Gamma_\infty^{\delta\#}.$ 
\qed

Denote the 1-form \eqref{HQE:descendents} corresponding to $\T=\D^M$
by $\Omega_\infty$ and the 1-form \eqref{HQE:ancestors} corresponding 
to $\T=\A^M_\tau $ by $\Omega_\tau.$ 
\begin{lemma} \label{forms:t_inf}
Let $r\in \Z$, $\q'$ and $\q''$ be such that $\widehat w_\infty'-\widehat w_\infty''=r.$
Then up to factors independent of $\gl$ the 1-forms $\Omega_\infty(\q',\q'')$
and $\Omega_\tau([S_\tau \q']_+,[S_\tau\q'']_+)$ coincide.
\end{lemma}
\proof
We compute 
$\(\widehat S_\tau^{-1}\tensor \widehat S_\tau^{-1}\)\Omega_\tau$ in 
two different ways: by  using formula  \eqref{action:s} and by commuting $\(\widehat S_\tau^{-1}\tensor \widehat S_\tau^{-1}\)$
through the vertex operators of $\Omega_\tau.$ 

In the first case, up to a factor independent of $\gl$  we get
$\Omega_\tau([S_\tau \q']_+,[S_\tau\q'']_+).$ 
In the second one, using Lemmas \ref{conjugation:a} and \ref{conjugation:d}
we get, up to factors independent of $\gl,$  the 1-form 
$\Omega_\infty(\q',\q'').$ \qed

{\em Proof of \thref{t1}.} Let $r\in \Z$ be arbitrary and assume that 
$\q',\q''\in \H_+$ are such that 
$\widehat w_\infty' - \widehat w_\infty'' =r,$ i.e.,
\ben
r\ge = \Omega (w_\infty',\q')- \Omega(w_\infty'',\q''). 
\een   
Since $S_\tau$ is a symplectic transformation, the RHS equals 
\ben
\Omega (w_\tau',S_\tau\q')- \Omega(w_\tau'',S_\tau\q'') = 
\Omega (w_\tau',[S_\tau\q']_+)- \Omega(w_\tau'',[S_\tau\q'']_+),
\een 
where the truncation operation  $[\ ]_+$ does not change the value 
of the symplectic form because $w_\tau\in \H_-.$  Theorem \ref{t1}
follows from \leref{forms:t_inf} and \thref{t2}.
\qed


\end{document}